\documentclass[12pt,reqno]{article}
\usepackage{graphics}
\newcommand{\ox}{\omega_x}
\newcommand{\oz}{\omega_z}
\newcommand{\Om}{\Omega}

\newcommand{\Omp}{\Omega^{per}_r}

\newcommand{\Ombp}{\Omega^{BP}_r}
\newcommand{\Ombpz}{\Omega^{BP}_0}
\newcommand{\GGbp}{{\cal G}^{BP}_r}
\newcommand{\GGbpz}{{\cal G}^{BP}_0}
\newcommand{\GGp}{{\cal G}^{per}_r}

\renewcommand{\div}{{\rm div\, }}
\newcommand{\curl}{{\rm curl\, }}

\newcommand{\RR}{{\bf R}}
\newcommand{\ZZ}{{\bf Z}}
\newcommand{\NN}{{\bf N}}

\newcommand{\eps}{\epsilon}
\newcommand{\PP}{{\cal P}_*}
\newcommand{\BP}{{\cal BP}(N,s,q,\vec k)}
\newcommand{\tH}{t{\cal H}}
\newcommand{\tHH}{\overline{t{\cal H}}}
\newcommand{\BPP}{{\cal BP}_*}
\newcommand{\BPPP}{{\cal BP}_*(N,s,q_m, m{\bf Z})}
\newcommand{\TT}{{\cal T}}
\newcommand{\SSS}{{\cal S}}
\newcommand{\Pf}{\noindent {\it Proof: \  }}
\newcommand{\QED}{\newline $\diamondsuit$}
\newcommand{\iint}{\int\!\!\!\!\int}
\newcommand{\meanh}{\langle h\rangle}

\newtheorem{thm}{Theorem}[section]
\newtheorem{lem}[thm]{Lemma}

\newtheorem{prop}[thm]{Proposition}
\newtheorem{rem}[thm]{Remark}

\newcommand{\be}{\begin{equation}}
\newcommand{\ee}{\end{equation}}
\newcommand{\bea}{\begin{eqnarray}}
\newcommand{\eea}{\end{eqnarray}}
\newcommand{\beann}{\begin{eqnarray*}}
\newcommand{\eeann}{\end{eqnarray*}}
\newcommand{\nnn}{\nonumber}

\begin{document}
%\begin{titlepage}
\title{Periodic vortex lattices for the Lawrence--Doniach \\ model
of layered superconductors \\ in
a parallel field}
\author{{\Large S. Alama\footnote{Dept. of Mathematics and Statistics,
McMaster Univ., Hamilton, Ontario, Canada L8S 4K1.  Supported
by an NSERC Research Grant.}, \ A.J. Berlinsky%
\footnote{Dept. of Physics and Astronomy, McMaster Univ., Hamilton,
Ontario, Canada L8S 4K1.  Supported by an NSERC Research Grant.}, \
\&  L. Bronsard${}^*$}}
\thispagestyle{empty}
\maketitle

\begin{abstract}
We consider the Lawrence--Doniach model for layered superconductors,
in which stacks of parallel superconducting planes are coupled
via the Josephson  effect.  We assume
that the superconductor is placed in an external magnetic field
oriented  {\it parallel}  to the superconducting planes and study
periodic  lattice configurations in the limit as the Josephson
coupling parameter $r\to 0$.  This limit leads to the
``transparent state'' discussed in the physics literature, which
is observed in very anisotropic high-$T_c$ superconductors at
sufficiently high applied fields and below a critical temperature.  We use a
Lyapunov--Schmidt reduction to  prove that energy minimization
uniquely determines the geometry of the optimal vortex lattice:
a period-2 (in the layers) array proposed by  Bulaevski\u\i \,
\& Clem.  Finally, we discuss the apparent conflict with previous
results for finite-width samples, in which the minimizer in the
small coupling regime takes the form of ``vortex planes'' (introduced
by Theodorakis and Kuplevakhsky.)

\end{abstract}

\newpage

%\tableofcontents
%\end{titlepage}

\baselineskip=17pt

\section{Introduction}

In this paper we continue the analysis of  layered superconductors in a
parallel external magnetic field started in our previous paper
\cite{ABB1}.  In the previous paper we considered the case of
a superconducting sample of finite width; here we treat periodic
solutions in an infinitely wide superconductor.

The Lawrence--Doniach model \cite{LD} is a mesoscopic
Ginzburg--Landau type model for superconducting materials
with a planar layered structure, and was originally applied to
study organic superconductors and other superconducting
composites manufactured by deposing successive thin layers
of superconducting metal with interposing insulating films.
Interest in this model has
been spectacularly revived by the discovery of  high
temperature superconductors, since nearly all of these
materials exhibit a distinctly layered structure.
Indeed, a pure monocrystalline sample of cuprate high-$T_c$
material (such as Bi${}_2$Sr${}_2$CaCu${}_2$O${}_8$
(BSCCO), 
Tl${}_2$Ba${}_2$CaCu${}_2$O${}_8$ (TBCCO),
or to a lesser extent YBa${}_2$Cu${}_3$O${}_7$ (YBCO)) consists of
copper oxide superconducting planes stacked with
intervening insulating (or weakly superconducting) planes.

\medskip

\noindent {\bf The Lawrence--Doniach model.}\
We model the layered superconductor as an infinite stack
of superconducting planes, each parallel to the $xy$-plane
and with uniform separation distance $p$.  The planes
are thus described by those
$(x,y,z)$ with $(x,y)\in\RR^2$ and $z=z_n:=np$, $n\in\ZZ$. 
We impose an applied ``external'' magnetic field, of
constant magnitude $H$ and lying parallel to the planes,
along the $y$-direction,
$\vec H= H\hat y$.  We make the ansatz that the local
magnetic field inside the sample will be everywhere independent of
$y$ and point in the $y$-direction, 
$$    \vec h(x,y,z)=h(x,z)\, \hat y.   $$
The vector potential $\vec A$ may then be chosen to lie
in the $xz$-plane,
$$  \vec A(x,y,z)= A_x(x,z)\, \hat x \, + A_z(x,z)\, \hat z,
\qquad  \vec h = \curl \vec A = 
    \left(  {\partial A_x\over \partial z} -
         {\partial A_z\over \partial x} \right) \, \hat y.  $$

We study configurations which are bi-periodic, in the sense
that the physically observable quantities (the magnetic
field, the density of superconducting electrons,
and their currents) are  doubly periodic functions in the
$xz$-plane. More precisely, given $q>0$, $N\in\NN$,
and $s\in\RR$ we define the fundamental domain of
periodicity to be the parallelogram $\Pi=\Pi_{N,s,q}$ spanned by
the vectors  $\vec e_1=(2q,0)$ and $\vec e_2=(s,Np)$:
$$  \Pi:= \{(x,z)=t_1e_1 + t_2 e_2, \ 0\le t_1,t_2\le 1\}.  $$
In other words, our fields and currents will be $2q$-periodic
in $x$ (along each superconducting plane), and will repeat
themselves after each $N$ planes but with a horizontal translation
of $s$.
We denote 
$$  z_n=np \quad \mbox{and} \quad x_n=ns/N, \quad n\in {\bf Z},  $$
so the part of the $n$th SC plane which lies within the basic
period module $\Pi$ is described by the points
$(x,z_n)$ with $x_n\le x\le x_n+2q$.

In each plane we define a (complex-valued)
superconducting order parameter
$\psi_n(x)$, $n\in\ZZ$.  We choose units in which
$|\psi_n|=1$ represents a purely superconducting
state.  With the assumption that the physically
observable quantities are $\Pi$-periodic we
may measure the free energy over a single period
to obtain the following functional,
\bea\label{gibbs}
\GGbp (\psi_n,\vec A)&=& {H_c^2\over 4\pi} \left\{
   p\sum_{n=1}^N
    \int_{x_n}^{x_n+2q}  \left[
        {1\over \kappa^2}
           \left|\left( {d\over dx} -iA_x\right)\psi_n  \right|^2 +
       \frac12 (|\psi_n|^2-1)^2  \right]\, dx 
         \right. \\   
 \nnn      &&\qquad + {r\over 2}\, p\sum_{n=1}^N
                  \int_{x_n}^{x_n+2q}  \left| \psi_n - 
                     \psi_{n-1}\exp\left(i\int_{z_{n-1}}^{z_n} A_z(x,s)\, ds\right)
                          \right|^2
                  \, dx \\  
\nnn       &&  \left.
            \qquad   +\ {1\over  \kappa^2}
                \iint_\Pi \left(  {\partial A_x\over \partial z} 
                    - {\partial A_z\over \partial x}
      - H\right)^2 \, dx\, dz \right\}.
  \eea
Note that
$\GGbp$ is expressed in non-dimensional units, chosen
 such that the in-plane penetration
depth $\lambda_{ab}=1$, $\kappa=\lambda_{ab}/\xi_{ab}$
is the Ginzburg Landau parameter,  $r$ is the
{\it interlayer coupling parameter}
(or {\it Josephson coupling parameter},)
and the magnetic fields
are measured in units of $H_c/\kappa$,
where  $H_c$ is the thermodynamic critical field. (See \cite{Tinkham}.)

The coupling between the superconducting planes given
by the second sum in $\GGp$ simulates the Josephson
effect, by which superconducting electrons travel
from one superconducting region to another by quantum mechanical
tunnelling.  We will see this explicitly in the Euler--Lagrange
equations, where the currents in the gaps between planes will
be determined by the sine of the gauge-invariant phase difference.
The interlayer coupling parameter $r$
 gives the strength of the Josephson coupling.

\smallskip

As in our previous study of the boundary-value problem for
the Lawrence--Doniach system \cite{ABB1} we will
study  the minimizers (and low-energy
solutions)  for $r$ near zero.
Indeed, for the high-$T_c$ cuprates at a temperature sufficiently
below $T_c$ we expect $r$ is a small parameter
in the Lawrence--Doniach model.  In particular, for BSCCO or
TBCCO we expect $r(0)\simeq 10^{-4}$--$10^{-3}$.  (See
``Physical background'' below.)

\medskip

\noindent {\bf  Results.}
We find that for $r\sim 0$ there is an {\it unique}
periodic solution of the Lawrence--Doniach system
which attains the minimum energy (per unit area)
among any other periodic configuration, with any
period geometry.  This absolute minimizing
solution has fundamental
domain $\Pi_{N,s,q}$ with $N=1$, $s=q=\pi/Hp$: in other words, 
the currents and field
are $2\pi/Hp$-periodic in $x$, and shifting from one plane to
the next results only in a horizontal translation
by a half-period $s=\pi/Hp$.  For example,
the magnetic field satisfies
\be\label{1periodic}
h\left(x+{2\pi\over Hp},z\right)=h(x,z)=h\left(x+{\pi\over Hp}, z+p\right).
\ee
Clearly this condition implies that the magnetic field, superconducting
electron density and currents are $2p$-periodic in $z$.
In addition, the flux per period for this solution is
exactly one quantum fluxoid, $2\pi$.
This configuration is optimal in the following sense:
each choice of period geometry and quantized flux
determines a function space in which to minimize
$\GGbp$.  The configuration described above
attains the minimum of free energy per unit
cross-sectional area among all possible choices
of $\Pi_{N,s,q}$ and flux quantization.  In fact our
result says more:  for any given geometry and flux,
then for all $r$ sufficiently small
the minimum energy per unit cross-sectional area 
in that class is either always
strictly larger than the $\Pi_{1,{\pi\over Hp},{\pi\over Hp}}$--bi-periodic
solution above, or exactly equal. The  energies
 coincide if and only if the given lattice $\Pi_{N,s,q}$ is
commensurate with $\Pi_{1,{\pi\over Hp},{\pi\over Hp}}$,
and its minimizer coincides with the  (unique)
$\Pi_{1,{\pi\over Hp},{\pi\over Hp}}$-minimizer described above.
  (See Theorem~\ref{bigthm}.)

The minimization among configurations and among geometries
and fluxes is entirely the result of direct rigorous analysis of the
Lawrence--Doniach system, without recourse to numerical approximation
or further restrictions other than the periodic ansatz.
This is not like the case for the Abrikosov lattice in the Ginzburg--Landau
model near $H_{c2}$, where the energy-minimizing geometry
and flux quantization were determined numerically by comparisons
within a finite collection of configurations.
Our result confirms the prediction of
Bulaevski\u \i  \ \& Clem \cite{BC}, who
claimed that for sufficiently large applied fields $H$
the energy minimizers should
form a period-two vortex lattice satisfying an
ansatz of the form (\ref{1periodic}).

\medskip

We also consider the case of finitely many superconducting
planes, with periodic observable quantities.  The result we
obtain is much the same as in the bi-periodic case:  minimization
of the free energy per period strip dictates the choice of
period $2q={2\pi\over Hp}$, and the minimizing solutions
are (to order $r$) $2p$-periodic in $z$ in the interior of
the sample.  (See Theorem~\ref{Thm2}.)

\smallskip

In both settings we obtain a ``transparent state'' as observed
in experiments on BSCCO by Kes, Aarts, Vinokur, \&
ver der Beek \cite{Kes}  (see ``Physical background'' below.)
  Superconductivity is essentially unaffected
by the presence of strong magnetic fields, $|\psi_n|= 1-O(r)$, while
the magnetic field is virtually unscreened by
the presence of the superconducting planes, $h=H+O(r)$.
In particular, the order parameters are never zero, and the
``vortices'', which correspond to the local maxima of the local
magnetic field, fit entirely in the gaps between the superconducting layers.
In the physics literature these are referred to as {\em Josephson
vortices} as opposed to the Abrikosov vortices characterized
by the vanishing order parameter in the Ginzburg--Landau model.

\medskip

The juxtaposition of the results obtained for the finite-width
samples in our previous paper \cite{ABB1} with the result stated
above may seem contradictory:  in \cite{ABB1} we prove that
for {\it any} sample of finite width $-L\le x\le L$ (and for $r\sim 0$)
the unique absolute energy minimizer occurs when the Josephson
vortices are aligned vertically in a ``vortex planes'' geometry.
These solutions were described by the physicists Theorodakis \cite{Th}
and Kuplevakhsky \cite{K}, the latter claiming that they were the
only solutions to the Lawrence--Doniach system.
Our explanation for this duality rests on the dependence of
the interval of validity of the perturbative regime $r\sim 0$
on the length of the sample size:  the larger the sample, the
smaller the interval in $r$ for which the results of \cite{ABB1}
are valid.  We will provide a more complete comparison
of the two solutions in section~8.

\bigskip

\noindent {\bf Physical background.}\
A first attempt to model layered superconductors 
superconductors is by an {\it anisotropic Ginzburg--Landau}
model, which treats the sample as a three-dimensional
solid with anisotropic material parameters.  In particular,
the  coherence length
along the perpendicular to the planes $\xi_c$
is assumed
to be different from the value of
the coherence length $\xi_{ab}$
within the planes (see \cite{Tinkham} or \cite{CDG}, for
example) and the anisotropy is measured by an
``effective mass ratio''  $\Gamma:= \xi_{ab}^2/\xi_c^2$.
For example Iye \cite{Iye} cites values of
 $\Gamma_{YBCO}\simeq 49$,  $\Gamma_{BSCCO}\simeq 3025$,
and  $\Gamma_{TBCCO}>90000$.
For certain materials (such as YBCO) and temperatures close
to the critical temperature $T_c$ this approximation
seems valid, but for more anisotropic superconductors
the anisotropic Ginzburg--Landau model does not give a
good qualitative or quantitative description of experimental
observations.   

Experimental evidence of the failure of the anisotropic
Ginzburg--Landau model was observed by Kes, Aarts,  Vinokur, and
 van der Beek \cite{Kes},  in the situation where the sample is placed in a
strong magnetic field oriented {\it parallel}\, to the superconducting
planes.  Experiments  reveal a crossover between ``three-dimensional''
behavior  (governed by the anisotropic Ginzburg--Landau model)
and ``two-dimensional'' behavior when the temperature is lowered
beyond a critical value  $T_{c0}< T_c$.
The ``two-dimensional'' regime is characterized by
a ``magnetically transparent'' state, in which the  applied magnetic
field penetrates completely
between the planes, virtually unscreened by the superconductor.
These observations are  not in agreement with the anisotropic
Ginzburg--Landau model, where the magnetic field is 
largely expelled from the bulk and isolated vortices appear
in a triangular ``Abrikosov lattice.''
For highly anisotropic materials (such as BSCCO and TBCCO)
the two-dimensional regime occurs within one degree Kelvin
of $T_c$, and an analysis of the transparent state requires a model which
addresses the discrete nature of the material: the 
Lawrence--Doniach model.  

\smallskip

To motivate the treatment of  $r$ as a small parameter, we
re-write $r$ in the original dimensional coordinates,
$$   r=r(T)={2\xi_{ab}^2 \lambda_{ab}^2\over
           \bar\lambda_J^2 \bar p^2}
           = {2\over \lambda_J^2 \kappa^2 p^2},
$$
 where
$\bar\lambda_J=\lambda_J\,\lambda_{ab}$, $\bar p = p\, \lambda_{ab}$
are the actual physical Josephson penetration depth and separation
distance.   The length scale
$\lambda_J$ in the Lawrence--Doniach model is directly related to the
effective mass ratio of the anisotropic Ginzburg--Landau model via
$\Gamma=\lambda_J^2$.
The temperature dependence of $r(T)$ is determined by
the conventional
dependences of $\xi^2_{ab}\simeq \xi^2_{ab}(0){T_c\over T_c-T}$ and
$\lambda^2_{ab}\simeq \lambda^2_{ab}(0){T_c\over T_c-T}$.
For the high-$T_c$ cuprates, 
$\xi_{ab}(0)$ and $\bar p$ are of the same order of magnitude, and hence
for highly anisotropic superconductors at temperatures sufficiently
below $T_c$,
$r\sim \lambda_J^{-2}=\Gamma^{-1}$ is small.  In particular, for BSCCO or
TBCCO we expect $r(0)\simeq 10^{-4},10^{-3}$.

We note that Chapman, Du, \& Gunzburger \cite{CDG} have proven
that solutions of the Lawrence--Doniach model converge to
solutions of the anisotropic Ginzburg--Landau model (and
in particular the convergence of energy minimizers) under the
limit $p\to 0$ with $\kappa,\Gamma$ fixed.  This limit does not 
correspond to our ``two-dimensional'' regime, since it
would send $r\to\infty$ and is therefore far from the ``weak coupling''
of superconducting planes observed in \cite{Kes}. 
Indeed, condition~(1) in \cite{Kes}, which defines the dimensional cross-over
point $T=T_{c0}$, is equivalent to $r(T_{c0})\le 1$ in the
Lawrence--Doniach model.
 On the other hand, by fixing
$\lambda_{ab}=1$ in our units the limit $T\to T_c$ (from below)
sends  $p=\bar p/\lambda_{ab}(T)\to 0$.  Therefore it is not
surprising to recover  ``three-dimensional'' behavior for
$p\to 0$.

\bigskip

\noindent {\bf Methods.}\
As in our previous study of the boundary-value problem for
the Lawrence--Doniach system \cite{ABB1} we will
use a degenerate perturbation approach to
study  the minimizers (and low-energy
solutions)  $r$ near zero.  However, the periodic and bi-periodic
settings are more subtle and present some additional complications
in applying the method.

First, periodic magnetic fields and currents  are generally represented by
non-periodic order parameters and potentials $(\phi_n,\vec A)$.
In order to define a variational setting for the bi-periodic problem
we adopt the idea of `t~Hooft \cite{tH} to define spaces of
functions which are periodic up to gauge transformation.
(See section~2.)  These spaces are rather complicated from
the functional analytic point of view, and therefore it is important
to find an equivalent formulation.  We do this by making a gauge
change (in section~3) to arrive at a family of affine Hilbert spaces
representing the periodic configurations.

\smallskip

For the finite-width case considered
in \cite{ABB1} the crucial observation was that when $r=0$ 
the planes decoupled, and the energy could be
minimized explicitly.   Even after gauge symmetries
were removed the  $r=0$ problem exhibited an additional
symmetry, corresponding to an $(N-1)$-dimensional torus action (where
$N$ denotes the number of superconducting planes,)  
and the minimization problem at $r=0$ degenerated on
 a finite dimensional hyperplane in function space.
When $r\neq 0$ this symmetry was broken and 
by a Lyapunov--Schmidt decomposition we reduced
the problem of finding solutions with $r\simeq 0$ to a
finite dimensional variational problem on this hyperplane.
In this finite-width setting, the minimum value of energy was
$O(r)$.

In the bi-periodic case the $r=0$ problem already
dictates a choice of period and flux:  the minimum energy
for the configuration $\Pi_{N,s,q}$ will be $O(r)$ if and
only if $q={m\pi\over Hp}$, $m\in\NN$, and the flux through $\Pi_{N,s,q}$ is
exactly $2\pi mN$.  (See Lemma~\ref{inf_per}.)
Within this more restricted class of
geometries  the
$r=0$ problem can again be solved explicitly, and the
same degenerate perturbation method as in \cite{ABB1}
applies to calculate the minimum energy as an expansion
in $r$ for $0<r<<1$.  Unlike the finite-width case, this
expansion yields a constant term at order $r$, and the computation
must be carried out to order $r^2$ in order to select
a point on the degenerate hyperplane which minimizes energy
with $r\sim 0$.
 The resulting finite dimensional minimization
problem may be completely resolved to give a clear choice
among the remaining parameters.  (See Theorem~\ref{Thm3}.)

\smallskip

The case of finitely many layers is similar (and simpler)
than the bi-periodic case;  we present a sketch in section~7
of the modifications necessary to treat that setting.

\medskip

As in \cite{ABB1} we recognize that in a real superconductor
the parameter $r$, while small, is not infinitessimal, which
raises the question of the validity of asymptotic results
in the regime $r\sim 0$ for physically relevent parameter
choices.  In section~5 of \cite{ABB1} we derived a lower
bound on the radius of validity of the $r$-expansion in
terms of the other parameters in the problem.  The same
analysis applies in the periodic problem, and we conclude
that the radius of validity is enhanced with smaller
period $q$ and $\kappa$, and larger external field $H$.
Since our optimal configuration ties the period $q$
to the reciprocal of the field $H$, we conclude that the
result is most applicable in strong external fields.
In particular, the expansion of the periodic problem
is more likely to be applicable than the expansion obtained
in \cite{ABB1}, since the validity of finite-width expansions 
deteriorates for large samples, whereas treating the
sample as infinitely large shifts the dependence to
the period, which can be small with large $H$.

\bigskip

We briefly outline the organization of the paper:
in the second section  we introduce a functional analytic
framework for the doubly periodic problem, based on
the elegant formulation of gauge periodicity due to `t~Hooft
\cite{tH}.  The third section defines an equivalent space, eliminating
bothersome symmetries due to gauge and continuous translation
invariance.  In section~4 we show how minimization at $r=0$
forces the correct choice of period and flux, and the fifth
section reviews the perturbation method used in \cite{ABB1}.
For the periodic problems it is necessary to expand the energy
to order $r^2$, and hence the solutions to the projected equations
(coming from the Lyapunov--Schmidt reduction) must be
calculated to order $r$:  these computations occupy
section~6.  The seventh section sketches the procedure for
the case of finitely many planes, which is conceptually simpler.

\smallskip

\noindent
{\bf Acknowledgement.}\
SA and LB wish to thank F. H\'elein and the Ecole Normale
Superieure--Cachan for their kind hospitality during 
the Winter and Spring of 2000.

\section{Variational formulation.}
\setcounter{thm}{0}

We begin with the doubly periodic case:  we adopt
the notation of the introduction, and assume that
the gauge invariant quantities (fields and currents)
are periodic in each SC plane, and allow the pattern of fields and currents
could repeat themselves with a horizontal shift in $x$ after
each period in $n$.  As before, for given $N\in\NN$, 
$s\in \RR$ and $q>0$ we define the fundamental domain
 $\Pi=\Pi_{N,s,q}$ as in the introduction.  
We seek solutions to the LD system which are periodic
with respect to the lattice generated by $\vec e_1,\vec e_2$,
\be
\label{hper}
h(x+2q,z)=h(x,z)=h(x+s,z+Np), 
\ee
and similarly for the other gauge-invariant quantities: the
density of superconducting electrons
$|\psi_n|$, the in-plane current density $j_x^{(n)}$, and
the Josephson current density in the gaps, $j_z^{(n)}$, defined
by
$$   j_x^{(n)} := 
    \mbox{Im}\, \left[
     \psi_n^*  \left( {d\over dx}-iA_x(x,z_n)\right)\psi_n  \right],
     \quad
     j_z^{(n)} :={r\kappa^2 p\over 2} \mbox{Im}\, 
         \left[\psi_n^*   
             \left(\psi_n - \psi_{n-1}e^{i\int_{z_{n-1}}^{z_n} A_z(x,z)\, dz}
                  \right)\right].
$$
Of course when $s=0$ we would have simple periodicity
in $z$, with period $Np$.

\bigskip

The subtlety of the periodic problems is that periodic
magnetic fields and currents
are generally represented by
{\it non-periodic} potentials $\vec A$ and order parameter $\psi_n$.
One setting for such periodic problems is via
 {\it t'Hooft boundary conditions} \cite{tH}, for which one demands that
 $\vec A$ and $\psi_n$ be periodic up to a family of gauge transformations
 from one period cell to the next.  
Let $N,q,s$ be fixed constants which
determine the period lattice $\vec e_1,\vec e_2$ as above.
We say that
 \be
 \label{tH}
 (\psi_n,\vec A)\in \tHH=\tHH (N,s,q)
 \ee
 if $\psi_n\in  H_{loc}^1(\RR)$, $\vec A \in H_{loc}^1(\RR^2;\RR^2)$ and
 there exist two
functions 
 $\omega_x$ and $\omega_z$ in
$H_{loc}^2(\RR^2)$ such that:
 \bea  \label{tHooft1}
 &&\psi_n(x+2q)=\psi_n(x)\exp[i\ox(x,z_n)] ,
 \quad \psi_{n+N}(x+s)=\psi_n(x)\exp[i\oz(x,z_n)],
 \\ && \nnn  \\  \label{tHooft2}
&&\vec A(x+2q,z)=\vec A(x,z)+\nabla\ox(x,z),\quad \vec A(x+s,z+Np)=\vec
A(x,z)+\nabla\oz(x,z).
 \eea

Since $\psi_{n+N}(x+2q)$ can be related back to $\psi_n(x)$
in two different ways (by commuting the order in which the two
rules in (\ref{tHooft1}) are applied,) this imposes a constraint
on the functions $\ox,\oz$ which are permissible in this
definition.  For any configuration $(\psi_n,\vec A)\in \tHH$
there exists an integer $K\in\ZZ$ such that
\be\label{compatibility}
  \ox(x,z_n)-\ox(x+s,z_{n+N})+\oz(x+2q,z_n)-\oz(x,z_n)
     =2\pi K.
\ee
Clearly, by continuity $K$ may be chosen independent of $x$;
we will see below that it is also independent of $n$.
The important role played by $K$ is revealed by
integrating the magnetic field strength over the basic unit
period cell, applying Stokes' Theorem, and substituting the
 identity (\ref{compatibility}):
\bea\label{flux1}
\iint_\Pi \curl \vec A \, dx \, dz 
   &=&  \int_{\partial\Pi} \vec A\cdot d\vec s  \\
        \nnn
        & =&  \ox(0,0)-\ox(s,Np)+\oz(2q,0)-\oz(0,0) \\
        \nnn
        &=&  2\pi K.
\eea
Hence the total flux per period cell is {\it quantized} for any 
element of $\tHH$.  Note also that the periodicity of
$h=\curl \vec A$ implies that the constant $K$ in (\ref{compatibility})
may indeed be chosen independent of $n$.

Since every $(\psi,\vec A)\in \tHH$ is associated to an integer $K$
via a continuous map,  
$$  (\psi,\vec A)\to {1\over 2\pi}\iint_\Pi \curl\vec A\, dx\, dz,  $$
this discrete choice separates the space $\tHH$ into disjoint connected
components,
$$  \tHH (N,s,q) = \bigcup_{K\in\ZZ} \tH(N,s,q,K).
$$
We may of course minimize the Lawrence--Doniach energy
in each component $\tH(N,s,q,K)$ separately, and indeed
 Lemma~\ref{inf_per} below will indicate the optimal choice
of $K$, for appropriate $(N,s,q)$.

\medskip

\begin{rem}\label{bundle}\rm
These periodic configurations may also be described within the
context of complex line bundles over the 2-torus ${\bf T}$, defined 
as $\RR^2$ modulo the discrete lattice generated by
$\vec e_1,\vec e_2$.  A configuration
$(\psi_n,\vec A)$ defines a connection on ${\bf T}$ by interpreting
$\vec A$ as a one-form, while the $\psi_n(x)$ are restrictions to the
$N$ planes in $\Pi$ of a section of a complex line bundle over
${\bf T}$, with structure group $U(1)$.  The `t~Hooft conditions incorporate the
$U(1)$ group action as we pass from one coordinate patch to another on the 
manifold ${\bf T}$.
 In this setting, the magnetic
field $h=\curl\vec A$ is the curvature of the connection, and 
(\ref{flux1}) is a Gauss-Bonnet relation between the total curvature
and the Euler number $K$, reflecting the fact that the bundle is nontrivial.
\end{rem}

\medskip

While this is a very elegant formulation of the periodic problem
it poses some practical problems for analysis of the variational
problem because of the introduction of the auxilliary functions
$\ox,\oz$ as well as the usual degeneracies associated
with gauge symmetry.  Fortunately we will be able to simplify
the setting of the problem in two ways.  First, by fixing a gauge
we find a much simpler Hilbert manifold setting.  This will
be done in section~3.  Second, we will show (see 
 Lemma~\ref{inf_per} below), that the least-energy solutions
 will have $|\psi_n(x)|\simeq 1$, and therefore it will
 be convenient to use  polar coordinates
for $\psi_n$ in order to treat the phases more directly.
In the remainder of this section we introduce subspaces
of $\tHH$ for $\psi_n$ in polar form and describe their
properties and the
Euler--Lagrange equations for critical points in these
spaces.

\medskip

 To define the space $\tHH$ in terms of $f_n$, $\phi_n$, 
with $\psi_ n(x)=f_n(x) e^{i\phi_n(x)}$ we must take into account that
 the periodicity conditions on $\phi_n$ can only hold modulo
 $2\pi$.  For the second relation in (\ref{tHooft1})
(translating in $z$ from $np\to (n+N)p$)
 this is not important, since the constant factor of $2\pi$ 
can be absorbed into $\oz$.  The degree
of winding per period in the $x$-direction plays a 
particularly important role,  and hence
 we assume that there exist integers $k_n$ such that
 \be\label{tHooft3} 
 \phi_n(x+2q) = \phi_n(x) + \ox(x,z_n) + 2\pi k_n,
    \quad \phi_{n+N}(x+s) = \phi_n(x) + \oz(x,z_n).  \ee
This definition does not uniquely determine
the winding numbers $k_n$, since adding the same constant
multiple of $2\pi$ to each merely reduces the value
of the function $\ox(x,z)$ by that same quantity.
To remedy this problem we may (without loss of generality)
fix the value
$$    k_0=0.  $$
Returning to the calculation (\ref{compatibility}), we see that
$$  K = k_{n+N}-k_n = k_N- k_0 =k_N,  $$
and hence the index $K$ associated to the space $\tH(N,s,q,K)$
measures the net change in the winding number of $\phi_n$
over one period $Np$ in $z$.
Of course the moduli of the order parameters is gauge-invariant, and
is $\Pi$-periodic,
\be\label{fper}  f_n(x+2q)=f_n(x)=f_{n+N}(x+s). 
\ee

\smallskip

We may now play the same game with the winding numbers
$k_n$ ($n\in\ZZ$) in the polar representation $(f_n,\phi_n,\vec A)$
 as we did with the single index $K$ for the configurations
$(\psi_n,\vec A)\in\tHH$.  Since $(f_n,\phi_n,\vec A)$ determines
$(k_n)_{n\in \ZZ}$  via a continuous selection,
 the space of
admissible configurations splits into disconnected classes.  
Fixing $\vec k:=(k_n)_{n\in \ZZ}$ (with $k_0=0$ and $k_{n+N}-k_n=K$)
we define our basic
space of bi-periodic configurations $\BP$ by choosing the corresponding
connected component of $\tHH$:
\beann
   \BP &:= & \{([f_n,\phi_n],\vec A): \
      f_n,\phi_n\in H^1_{loc}(\RR), \ \vec A\in H^1_{loc}(\RR^2;\RR^2), 
\ \mbox{and there}\\
&& \qquad \mbox{exist $\ox,\oz\in H^2_{loc}(\RR^2)$
such that (\ref{fper}), (\ref{tHooft3}), and (\ref{tHooft2})
hold.} \}.
\eeann
Each component $\BP$ describes a subclass of $\tH(N,s,q,K)$
with $K=k_N$, and therefore the flux quantization
$$  \iint_\Pi \curl\vec A\, dx\, dz = 2\pi K=2\pi k_N  $$
holds for every configuration $(f_n,\phi_n,\vec A)\in BP$.
As for $\tH(N,s,q,K)$, it will turn out that there are preferred classes $\BP$
 for which the free energy will be smallest possible.
  (See Lemma~\ref{inf_per} below.)

\medskip

It is a simple calculation to verify that for any configuration
$(\psi_n,\vec A)\in \tH(N,s,q,K)$ (and for
$(f_n,\phi_n, \vec A)\in \BP$) we recover the desired
periodicity conditions (see (\ref{hper}))
for the gauge-invariant quantities.  
Since only the gauge-invariant quantities appear in 
an expression of the free energy density,  the energy density will be
 $\Pi$-periodic and hence
we are justified in measuring the energy of the configuration
over only one period cell.  We may therefore define the Gibbs free
energy $\GGbp$ for the configuration $(\psi_n,\vec A)\in \tH(N,s,q,K)$
as in the introduction.
When $\psi_n$ can be represented by
polar coordinates $f_n$, $\phi_n$ we use the equivalent form:
 \beann
\Ombp (f_n,\phi_n,\vec A)&=&  p\sum_{n=1}^N \int_{x_n}^{x_n+2q}\left[
       \frac 12  (f_n^2-1)^2 + {1\over\kappa^2} (f'_n)^2
            + {1\over\kappa^2}(\phi'_n- A_x(x_n,y))^2 f_n^2
             \right]\, dx \\   
             &&\qquad + {r\over 2}  p\sum_{n=1}^N
                \int_{x_n}^{x_n+2q} \left( f_n^2 +f_{n-1}^2 
              -2 f_n f_{n-1}\cos(\Phi_{n,n-1} )\right) \, dx \\  
            &&\qquad  +\ {1\over  \kappa^2}
                \iint_\Pi \left(  {\partial A_x\over \partial z} 
                    - {\partial A_z\over \partial x}
      - H\right)^2 \, dx\, dz ,
  \eeann 
where $\Phi_{n,n-1}(x)$ is the {\it gauge invariant phase difference},
$$\Phi_{n,n-1}(x):= \phi_n(x)-\phi_{n-1}(x) 
      -\int_{z_{n-1}}^{z_n} A_z(x,z)\, dz.  $$
We observe that $\Phi_{n,n-1}(x)$ is  periodic modulo
$2\pi$,
$$  \Phi_{n,n-1}(x+2q)= \Phi_{n,n-1}(x) + 2\pi(k_n-k_{n-1}).  $$

\begin{rem}\label{gabriella}\rm
Note that although we have
defined  $(f_n,\phi_n,\vec A)\in\BP$ globally
in $\RR^2$, in fact the
gauge-invariant  quantities $f_n$, $\Phi_{n,n-1}$, 
$V_n:=(\phi'_n-A(x,z_n)$, and $h(x,z)$
are globally determined by the values of
 $f_n(x)$, $\phi_n(x)$, and $\vec A(x,z)$
for $n=1,\dots,N$ and  $(x,z)\in \Pi$, and
the values of $\oz(x,0)$.  This is very clear
for $f_n$, $V_n$, and $h$; for $\Phi_{n,n-1}$, 
we require $\phi_0(x)$ to determine $\Phi_{1,0}(x)$, and
it is here that the `t~Hooft condition comes into play
via the function $\oz(x,0)$.
In retrospect, we could have defined $(f_n,\phi_n,\vec A)$ only
on the period module $\Pi$ and $\ox,\oz$ only in a neighborhood
of the  left and bottom edges of $\Pi$.
This would be enough to make sense of the free energy, 
 and it allows interpretation of the
't~Hooft conditions as true boundary conditions,
\beann
 &&\psi_n(x_n+2q)=\psi_n(x_n)\exp[i\ox(x_n,z_n)] ,
 \quad n=1,\dots,N;
  \\  
&&
\psi_{N}(x+s)=\psi_0(x)\exp[i\oz(x,0)], \quad 0<x<2q;  \\
&& 
\vec A(x+2q,z)=\vec A(x,z)+\nabla\ox(x,z),\quad 
   0<x<s, \  z={Np\over s}x; \\
&&
\vec A(x+s,Np)=\vec
A(x,0)+\nabla\oz(x,0), \quad 0<x<2q.
 \eeann
 This is the interpretation of the 't Hooft conditions taken by
 Tarantello \cite{Ta} for example. \QED
\end{rem}

\bigskip

The Euler--Lagrange equations under 't Hooft conditions
are virtually identical to those derived in \cite{ABB1}
except that we must take into account
the  periodic boundary conditions 
 both in $x$ and $z$.  To obtain the equations
we may choose smooth genuinely $\Pi$-periodic test functions
to vary each unknown in turn.
We obtain:
\bea\label{bfeqn}
 &&-{1\over \kappa^2} f''_n + (f_n^2-1)f_n 
    + {1\over\kappa^2}(\phi'_n-A_x(x,z_n)^2 f_n\\
\nnn &&\qquad = {r\over 2}\left(
       f_{n-1}\cos\Phi_{n,n-1}
      +f_{n+1}\cos\Phi_{n+1,n}-2f_n\right); \\
\label{bhz}
&&{\partial h\over \partial z}(x,z)=0 \qquad z\neq z_n;\\  
\label{bjump}
&& h(x,z_n+)-h(x,z_n-) = - p f_n^2(x) (\phi'_n - A_x(x,z_n));\\
\label{bhx}
&&{\partial h\over \partial x} =
     {r\kappa^2 p \over 2} f_n(x) f_{n-1}(x) \sin\Phi_{n,n-1}(x),
            \mbox{ for $z_{n-1}<z<z_n$}.
 \eea
     These equations should be solved together with $\Pi_{N,s,q}$-periodic
boundary conditions for the observable quantities, $f_n(x)$, $h(x,z)$, 
$V_n(x):=(\phi'_n(x) - A_x(x,z_n))$, and $\sin\Phi_{n,n-1}(x)$.
  For later purposes, we also record the current conservation
  equation
  \be
  \label{cu}
  {1\over \kappa^2}{d\over dx}\left( f_n^2(\phi'_n - A_x(x,z_n))\right) =
   {r\over 2}(f_n  f_{n-1}\sin\Phi_{n,n-1}
       -f_{n+1}f_n \sin\Phi_{n+1,n}).
       \ee
This equation is not independent of the others:  it can be derived
by differentiating (\ref{bjump}) and substituting from (\ref{bhx}).
This is not surprising since $\phi_n$ and $\vec A$ are related
through gauge invariance.  We note that 
$$ j_x^{(n)}(x)=f_n^2(x)(\phi'_n(x)-A_x(x,z_n))  $$
measures the current density within the $n^{th}$ superconducting
plane, while 
$$  j_z^{(n)}(x)= {r\kappa^2 p\over 2} f_n(x) f_{n-1}(x)
               \sin\Phi_{n,n-1}(x)  $$
gives the Josephson current density in the gap between the
$(n-1)^{st}$ and $n^{th}$ planes.  In this way we may view
(\ref{cu}) as a semi-discrete version of the classical
equation of continuity $\div \vec j=0$.  This is the conservation
law corresponding to the $U(1)$ gauge invariance as guaranteed
by Noether's Theorem.

We also record a useful formula for $\Phi_{n,n-1}$ which follows
from Stokes' Theorem:
\be\label{Phieqn}
\Phi_{n,n-1}(x) = \int_0^x  (V_n-V_{n-1})\, d\bar x
   + p \int_0^x  h^{(n)}(\bar x)\, d\bar x + \Phi_{n,n-1}(0).
\ee

\section{Fixing a gauge}
\setcounter{thm}{0}

While the spaces $\tH (N,s,q,K)$ and $\BP$
 are indeed unwieldy for the purposes
of analysis, we will show that by fixing an appropriate gauge 
we can, without loss of generality, work in a much simpler setting.
We begin by noting that the periodic problem has a larger symmetry
group than the fixed interval problem studied in the previous
 paper \cite{ABB1}.  In addition to electromagnetic gauge invariance,
$$  f_n\to f_n, \quad \phi_n(x)\to \phi_n(x)-\lambda(x,z_n), \quad
   \vec A\to \vec A -\nabla\lambda(x,z),  $$
there is also {\it translation invariance}, 
$$  (f_n,\phi_n,\vec A)\to 
   (f_n(\cdot-x_0), \phi_n(\cdot-x_0),\vec A(\cdot-x_0,\cdot)),
   \qquad x_0\in {\bf R}.  $$
While this last symmetry is only one-dimensional, for our
degenerate perturbation theory it is convenient to eliminate all but
the essential symmetries of the $r=0$ problem which are broken
when $r\neq 0$, and so our choice of ``gauge'' will fix a
translation as well.
Note that there is also a discrete translation invariance in
$z$, $z\to z- kp$, $k\in \ZZ$, but this symmetry
(being discrete) will not create analytical difficulties for
our method and hence plays a less important role.

First consider the spaces $\BP$ for which $\vec k\neq \vec 0$.
(Recall that we have fixed $k_0=0$ in the definition of $\BP$.)
By the discrete translation invariance
in $z$ we may relabel the $z$-axis if necessary in order to
obtain
\be\label{k_1}
k_1\neq 0 = k_0. 
\ee
We observe that the average
value of the local magnetic field is fixed by the choice of
the space $\BP$:  indeed by (\ref{flux1}),
\be\label{flux2}  \meanh:= {1\over 2qNp}\iint_\Pi h(x,z) \, dx\, dz 
    = {\pi\over pq}{K\over N}\ ,
    \ee
where $h=\curl \vec A$.

Next we define an appropriate space which eliminates gauge
and translation symmetries.  Assume again that $\vec k\neq \vec 0.$
We say that $(f_n,\phi_n,\vec A)\in \BPP(N,s,q,\vec k)$ 
provided there are constants
$\omega,d\in\RR$
such that
\begin{eqnarray}
 \label{BP1}
&& \vec A\in H^1_{loc}(\RR^2;\RR^2), \  f_n\in H^1_{per}(\RR), \ 
        \phi_n\in H^1_{loc}(\RR); \\ &&\quad \nnn \\
  &&\label{phi}
    \phi_n(x+2q)=\phi_n(x) + \omega+ 2\pi k_n;
 \\ && \quad \nnn\\
\label{BP4}
&&
\phi_{n+N}(x+s)=\phi_n(x)+{k_N\pi\over q} x +d ; 
\\ && \quad \nnn\\
\label{BP5}
&& 
f_n(x+2q)=f_n(x)=f_{n+N}(x+s); \\
&& \nnn \\
\label{BP6}
&&  \left.\begin{array}{l}
  \vec A(x,z)=\meanh\, ( z,0) + 
\left(\partial_z \xi \ , \,
       -\partial_x \xi\right), \quad\mbox{where} \\
     \\
  \xi\in H^2_{loc}(\RR^2), \quad \xi(x+2q,z)=\xi(x+s,z+Np)=\xi(x,z) ;
\end{array} \right\}
\\ \nnn && \\
\label{P4}
 && \int_{-q}^q \phi_0(x)\, dx =0=\int_{-q}^q \phi_1(x)\, dx.
\end{eqnarray}
As usual, we may choose a definition of the phases so that
\be\label{normalize}
\mbox{$0\le \phi_n(0)<2\pi$ for each $n\in \ZZ$.}
\ee
In the case $\vec k=\vec 0$, we define ${\cal BP}_* (N,s,q,\vec 0)$
via (\ref{BP1})--(\ref{BP6}), and we demand only
the vanishing of the mean value of $\phi_0(x)$ in (\ref{P4})
(with no restriction on the mean of $\phi_1(x)$).

It is easy to see that 
$\BPP(N,s,q,\vec k)\subset\BP$.
We must show that the two spaces are equivalent, in the sense
that for any configuration in $\BP$  there is
a gauge transformation which sends it to an element of 
$\BPP(N,s,q,\vec k)$.   In particular, this implies that
minimization of $\Ombp$ in $\BPP(N,s,q,\vec k)$  is equivalent
to minimization of $\Ombp$ in $\BP$. 
In addition, we will show that the
choice of gauge in $\BPP(N,s,q,\vec k)$ is a ``Coulomb'' gauge, in
the sense that the  $H^1(\Pi)$ norm of $\vec A$ is controlled 
by its curl in $L^2(\Pi)$. 

 \begin{thm}\label{bipercoulomb}  Let $N,s,q,\vec k$ be fixed.

\noindent
(a)\  There exists a
constant $C_0>0$ such that
for any $(f_n,\phi_n,\vec A)\in\BPP(N,s,q,\vec k)$,
$$  \|\vec A\|^2_{H^1(\Pi)} \le C_0 \|h\|^2_{L^2(\Pi)}, $$
where $\Pi$ is the period module.

\noindent
(b)\
For any $(f_n,\phi_n,\vec A)\in \BP$ 
there exists $\lambda\in H^2_{loc}(\RR^2)$ and $x_0\in \RR$
such that 
$$
\left( f_n(x-x_0),  \phi_n(x-x_0) - \lambda(x-x_0,z_n), 
\vec A(x-x_0,z) - \nabla \lambda(x-x_0,z)\right)\in\BPP(N,s,q,\vec k).
$$
\end{thm}
When $\vec k=\vec 0$ the choice of $x_0$ in (b) is immaterial:
we cannot remove translation invariance in that case.

\noindent
{\it Proof of Theorem~\ref{bipercoulomb}:}\

Suppose $(f_n,\phi_n,\vec A)\in\BPP$, with
$\vec A$ defined as in (\ref{BP6}).  Then $\Delta\xi = h-\meanh$
in $\RR^2$, and therefore the estimate follows from standard
elliptic regularity theory.

\smallskip

To prove (b), assume $(f_n,\phi_n,\vec A)\in\BP$ and
let $\xi$ be a solution to the periodic problem
\be\label{elliptic}  \left. \begin{array}{l}
     \Delta\xi = \curl\vec A - \meanh, \\
     \xi(x+2q,z)=\xi(x,z) =\xi(x+s,z+Np),
\end{array} \right\}  \ee
and $\tilde A=(\partial_z\xi, -\partial_x\xi)$.
By standard elliptic regularity theory $\xi\in H^2_{loc}(\RR^2)$.  Then
$\curl (\tilde A-\vec A)=0$
so that there exists $\hat\lambda \in H^2_{loc}(\RR^2)$
with $\tilde A=\vec A -\nabla\hat\lambda$.  
In fact $\hat\lambda$ is only determined up to a constant; 
set $\lambda=\hat\lambda -c$, with $c$ as yet undetermined.

Define 
$$ \tilde\phi_n:=\phi_n-\lambda(x,z_n).  $$
It follows that
\be
\label{phitilde}
\tilde\phi_n(x+2q)-\tilde\phi_n(x)=
\ox(x,z_n)-\lambda(x+2q,z_n)+\lambda(x,z_n)+2\pi k_n.
\ee
Since $\tilde A$ is periodic in $x$, we have (for all $z$)
$$\nabla\ox(x,z)
-\nabla \lambda(x+2q,z)+\nabla\lambda(x,z)=0,$$
and hence there exists a constant of integration $\omega$ so that
$$\ox(x,z)-\lambda(x+2q,z)+\lambda(x,z)=\omega $$
for all $(x,z)\in \RR^2$.
Substituting this identity in (\ref{phitilde}) we verify (\ref{phi}) 
for $\tilde\phi_n$.
Similarly
$$\tilde\phi_{n+N}(x+s)-
\tilde\phi_n(x)=\oz(x,z_n) -\lambda(x+s,z_n+Np)+\lambda(x,z_n),$$
and since
$$\tilde A(x+s,z+Np)-\tilde A(x,z)=
(Np,0)\meanh $$
it follows that there is another constant of integration $d$ such that
$$\oz(x,z) -\lambda(x+s,z+Np)+\lambda(x,z)=
Np\meanh \, x + d $$
for all $(x,z)\in \RR^2$.
In particular using the quantization formula (\ref{flux2}), it
follows
$$\tilde\phi_{n+N}(x+s)-
\tilde\phi_n(x)={K\pi\over q} x+  d, $$
and (\ref{BP4}) is satisfied.

\smallskip

To complete the argument in case $\vec k\neq \vec 0$, we define
$$  g(t):= \int_{t-q}^{t+q} \left(\tilde\phi_1(x)-\tilde\phi_0(x) \right)\, dx.  $$
Recalling (\ref{k_1}), $k_0=0\neq k_1$, so (\ref{phi}) implies that
$\tilde\phi_1(x)-\tilde\phi_0(x)\sim {k_1\pi\over q}x$
for $|x|$ large.  Therefore $g(t)\to\pm\infty$ as $t\to\pm\infty$,
or vice-versa.  Since $g(t)$ is continuous there exists $x_0\in \RR$ so that
$g(x_0)=0$.  We now choose the constant term in $\lambda$ to be
$$  c={1\over 2q}\int_{x_0-q}^{x_0+q} \tilde\phi_1(x)\, dx =
    {1\over 2q}\int_{x_0-q}^{x_0+q} \tilde\phi_0(x)\, dx,  $$
 and the conclusion
of (b) follows with $\lambda$ the desired gauge change 
and $x_0$ the translation.

In case $\vec k=\vec 0$, we cannot assert the existence of
such a translation $x_0$, but we can still choose 
$c={1\over 2q}\int_{-q}^q \phi_0(x)\, dx$ to obtain 
 $\tilde\phi_0$ with zero mean in $[-q,q]$.
\QED

\medskip

The spaces $\BPP(N,s,q,\vec k)$ are quite reasonable from the
analytical point of view.  First, note that each $\phi_n(x)$ is
the superposition of a linear function (with slope $\omega + 2\pi k_n$) with
a $2q$-periodic function, and  hence $\phi'_n(x)$ is $2q$-periodic.
Since 
$$  \omega=\phi_1(x_1+2q)-\phi_1(x_1)-2\pi k_1
\quad
\mbox{and}\quad
   d={1\over 2q}\int_{s-q}^{s+q} \phi_N(x)\, dx, $$
by the conditions (\ref{phi}), (\ref{BP4}) the quantity
$$  \sum_{n=1}^N \|\phi_n\|^2_{H^1[x_n,x_n+2q]}  $$
controls the $H^1_{loc}(\RR)$ norms of any finite collection
of $\phi_n$, $n\in\ZZ$.

Indeed, for any choice of parameters 
$\BPP(N,s,q,\vec k)$ forms an affine Hilbert space.
The tangent space (of admissible variations to 
$(f_n,\phi_n,\vec A)\in \BPP$) is the space
  $E=E(N,s,q)$ consisting of all $(u_n,v_n,\vec a)$ such that
\begin{enumerate}
\item[(1)] 
$u_n,v_n\in H^1_{loc}(\RR)$ ($n\in \ZZ$) and there exist
constants $\omega, d\in \RR$ such that
$$   u_n(x+2q) =   u_n(x)  = u_{n+N}(x+s),  $$
$$
 v_n(x+2q)=v_n(x) + \omega, \quad   v_{n+N}(x+s)=v_n(x) + d.  $$
\item[(2)]
$\vec a=(\partial_z \xi, -\partial_x \xi)$
for some $\xi\in H^2_{loc}(\RR^2)$ with
$$    \xi(x+2q,z) = \xi(x,z) = \xi(x+s,z+Np).  $$
\end{enumerate}
In view of the above remarks on controlling the extension
of elements of $\BPP(N,s,q,\vec k)$ beyond the basic
period $\Pi$, we may choose  a norm 
on $E$ of the form
$$ \|(u_n,v_n,\vec a)\|_E^2 =
   \sum_{n=1}^N \left[ \|u_n\|_{H^1([x_n,x_n+2q])}^2
        + \|v_n\|_{H^1([x_n,x_n+2q])}^2\right]
          + \iint_\Pi [\curl \vec a]^2\, dx\, dz. $$

It is clear that
$\Ombp$ is a smooth ($C^\infty$) functional on 
$\BPP(N,s,q,\vec k)$.  Since any configuration in $\BP$ 
can be associated to an element in $\BPP(N,s,q,\vec k)$ (via gauge
transformation and translation) with the same energy,
the minimizer of $\Ombp$ in $\BPP(N,s,q,\vec k)$ will also minimize
$\Ombp$ in $\BP$, and hence satisfy the Lawrence--Doniach
system of equations with periodic conditions.

\medskip

Finally, we note that the same procedure may be used to
fix a gauge in the space $\tH(N,s,q,K)$ of doubly gauge-periodic
configurations $(\psi_n,\vec A)$, to obtain a nice space \linebreak
$\tH_*(N,s,q,K)$
which eliminates gauge invariance.   Indeed, we say 
$(\psi_n,\vec A)\in \tH_*(N,s,q,K)$ if: 
\begin{enumerate}
\item[(a)]
$\psi_n\in H^1_{loc}(\RR)$ and
$\vec A\in H^1_{loc}(\RR^2;\RR^2)$; 
\item[(b)] $\vec A$ satisfies (\ref{BP6}); 
\item[(c)]
there exist constants
$\omega,d\in\RR$ such that 
$$
\psi_n(x+2q)=\psi_n(x)\exp[i\omega], \quad
     \psi_{n+N}(x+s)=
           \psi_n(x)\exp\left[i\left({K\pi\over q} x +d\right)\right].  $$
\end{enumerate}
A version of  Lemma~\ref{bipercoulomb} 
holds true for spaces $\tH_*(N,s,q,K)$:  for any
$(\psi_n,\vec A)\in \tH(N,s,q,K)$ there exists $\lambda\in H^2_{loc}(\RR^2)$
such that
$$   \left( \psi_n\exp(i\lambda(\cdot,z_n)), \vec A - \nabla\lambda \right)
                    \in \tH_*(N,s,q,K). $$
And for any $(\psi_n,\vec A)\in \tH_*(N,s,q,K)$ the
estimate of Lemma~\ref{bipercoulomb}~(a) also holds.

\section{Minimization at $r=0$}
\setcounter{thm}{0}

We are ready to begin the process of solving for the
lowest-energy periodic solution to the Lawrence--Doniach
system.  As we have already seen the periodic problems
are indexed by several parameters:  the number of
layers in the period module $N$; the geometry
of the period parallelogram $\Pi$ represented by
the
horizontal period
along the SC layers $2q$ and the horizontal shift $s$ when
advancing $N$ planes in $z$;
and the winding numbers $k_n$, which we will denote
collectively by $\vec k$.  It will turn out that
the geometry of the vortex lattice will be completely determined
by energy minimization!

\medskip

Denote the minimum free energy per unit cross-sectional area
over {\it all} admissible lattice configurations in $\BPP (N,s,q,\vec k)$ by
$$  \eps_r(N,s,q,\vec k):=
   \inf\left\{
      {\Ombp(f_n,\phi_n,\vec A)\over 2qNp}
      : \ (f_n,\phi_n,\vec A)\in 
          \BPP (N,s,q,\vec k) \right\}.  $$ 
Our goal is to minimize twice:  first for each {\em fixed} period geometry
$\Pi_{N,s,q}$ and choice of winding numbers $\vec k$, and
then to find the parameters $N,s,q,\vec k$ which give the least
energy per unit cross sectional area.
The main result we will prove is that, for $0<r<<1$, the periodic
solution with the least free energy among all possible
lattice geometries is a period-$2p$ in $z$
lattice, with period $2q={2\pi\over Hp}$ in $x$.
We denote by $\ZZ$ the integer sequence,
$k_n=n$ in the following statement of our theorem:
\begin{thm}\label{bigthm}
For any choice of $(N,s,q,\vec k)$ there
exists $\hat r=\hat r(N,s,q,\vec k)$ such that
either: 
$$ \eps_r(N,s,q,\vec k)>\eps_r(1,{\pi\over Hp},{\pi\over Hp}, \ZZ) \quad
\mbox{for all $r<\hat r$};  $$
 or 
$$\eps_r(N,s,q,\vec k)=\eps_r\left(1,{\pi\over Hp},{\pi\over Hp}, \ZZ\right)
\quad \mbox{for all $r<\hat r$,}
$$
in which case the minimizers are achieved
in $\BPP(N,s,q,\vec k)$, and they coincide with the period-$2p$
in $z$, period $2q={2\pi\over Hp}$ in $x$  minimizers in
$\BPP(1,{\pi\over Hp},{\pi\over Hp}, \ZZ)$.
\end{thm}

To completely clarify the above statement, we note that the
magnetic field and the supercurrents $j_x^{(x)}$, $j_z^{(n)}$
associated to elements of
the space $\BPP(1,{\pi\over Hp},{\pi\over Hp}, \ZZ)$
make a {\it horizontal translation} of a half-period $q_1:={\pi\over Hp}$
when we increase $z$ by its period $p$:
\be\label{shiftper}
\left.\begin{array}{c}
 f_n(x+2q_1)=f_n(x)=f_{n+1}(x+q_1), \quad 
V_n(x+2q_1)=V_n(x)=V_{n+1}(x+q_1)
\\  \\
 h(x+2q_1,z)=h(x,z)=h(x+q_1,z+p), \\  \\
          \Phi_{n,n-1}(x+2q_1)-2\pi= \Phi_{n,n-1}(x)= \Phi_{n+1,n}(x+q_1),
\end{array}\right\}
 \ee
where we recall $V_n(x)=(\phi'_n(x)-A_x(x,z_n))$.
It is clear that these solutions are also $2p$-periodic
in $z$.  These solutions have the period structure proposed
by Bulaevski\u\i  \, \& Clem \cite{BC}.  In section~6 we will
show that the associated  local magnetic field $h(x,z)$, in-plane and 
Josephson currents $j^{(n)}_x, j^{(n)}_z$ and density of superconducting
electrons $f_n$ admit an expansion near $r=0$ of the
form:
\be\label{per_2}  \left.  \begin{array}{c}
 h(x,z)=H + r\ (-1)^{n+1} 
{\kappa^2\over 2H}\cos(Hpx)\  +\  O(r^2), \qquad
    z_{n-1}<z<z_n,\\
\\
j_x^{(n)}(x)=r\ (-1)^{n+1}
     {\kappa^2\over Hp}\cos(Hpx) 
\ + \  O(r^2),\\
\\
j_z^{(n)}(x)= {1\over 2}\,r\, (-1)^n\, \kappa^2 p\sin (Hpx) \ + \ O(r^2),
\\  
\\
  f_n(x) = 1 \, - \, {r\over 2}\,  + \, O(r^2).
\end{array} \right\}
\ee

\medskip

The proof of Theorem~4.1 will be accomplished in many steps,
 concluding at the end of section~6.  The remainder of this section
is devoted to the study of the $r=0$ problem.  Section~5 introduces
a degenerate perturbation method based on a Lyapunov--Schmidt
decomposition at $r=0$, and section~6 solves the reduced problems
arising from the decomposition.

\bigskip

\noindent
{\bf   Choosing the period and winding numbers.}\
The first step in proving Theorem~\ref{bigthm} is to
consider the case $r=0$, in which case the superconducting
planes decouple.  We know from the previous treatment of
the finite-width problem in \cite{ABB1} that the case
$r=0$ is analogous to the self-dual point $\kappa=1/\sqrt{2}$
in the Ginzburg--Landau model in the sense that 
at that point minimizers satisfy a first-order Bogomolnyi 
system, in addition to the (second-order) Euler--Lagrange equations.
Indeed, we will see that the reduction to a first-order system
of equations is only possible
when the period $q$ and the winding numbers $\vec k$ are chosen
appropriately.  For a generic choice of $q,\vec k$ the minimum
of energy will be much larger.

\begin{lem}\label{inf_per}   For $r=0$,
\begin{enumerate}
\item[(a)] $\inf\{ \GGbpz (\psi_n,\vec A): \
    (\psi_n,\vec A)\in \tH_*(N,s,q,K)\} =0$ if and
only if there exists $m_0\in {\bf N}$ such that $Hpq=m_0\pi$
and $K=m_0 N$.
\item[(b)]  Assume $Hpq=m_0\pi$ and $K=m_0 N$, $m_0\in {\bf N}$.
For any $r\ge 0$ we have
$$  \inf\{\GGbp (\psi_n,\vec A): \
    (\psi_n,\vec A)\in \tH_*(N,s,q,K)\} \le {H_c^2\over 2\pi}Nqp\, r.  $$
Moreover, there exist constants $r_0$, $C>0$ such
that whenever $\GGbp(\psi_n,\vec A)\le {H_c^2\over 2\pi}Nqp\, r$, 
and $0\le r<r_0$, then $|\psi_n(x)|\ge 1- Cr^{1/2}>0$ for all $x$, $n$.
\item[(c)]
$$  \inf_{\BPP (N,s,q,\vec k)} \Ombpz = 0
$$
if and only if there exists  $m_0\in {\bf N}$ such that
$Hpq= m_0\pi$ and $k_n=m_0 n$ for each $n\in \ZZ$.
\end{enumerate}
\end{lem}
We remark that the constants $C, r_0$ in (b) will depend
on $N,s,H,m_0,\kappa$.

\Pf
We begin with (a), (c).
Assume $q={m_0\pi\over Hp}$, $k_n=m_0 n$,
$K=k_N=m_0 N$, with $m_0\in {\bf N}$, 
and choose
\be\label{minzero} f_{n0}\equiv 1, \quad
    \phi_{n0}=\alpha_n + nHpx, \quad \psi_{n0}=e^{i\phi_{n0}},\quad
 \vec A_0(x,z)=(Hz,0), 
\ee
for any $\alpha_n\in\RR$.  It is easy to see that
 $(f_{n0},\phi_{n0},\vec A_0)\in \BPP (N,s,q,\vec k)$, 
$(\psi_{n0},\vec A_0)\in \tH_*(N,s,q,K)$, with
$\GGbpz(\psi_{n0},\vec A_0)=0$ and
$\Ombpz(f_{n0},\phi_{n0},\vec A_0)=0.$  

To see the opposite
implications, assume first that
$\inf_{\BPP (N,s,q,\vec k)} \Ombpz = 0.$
By a simple estimate (see Proposition~2.1 of \cite{ABB1})
any minimizing sequence must have $f_n\to 1$ uniformly.
It is then easy to see that the energy controls the norm of the
minimizing sequence in $\BPP$, and the minimizing
sequence converges strongly to 
$(\hat f_{n},\hat\phi_{n},\hat A)\in\BPP (N,s,q,\vec k)$,
and attains the minimum of zero energy.
Being a sum of positive terms, $\Ombpz(f_{n0},\phi_{n0},\vec A_0)=0$
implies that each term is individually zero,
$$   f_n=0, \quad \phi'_n-A_x(x,z_n)=0,
       \quad  \curl \vec A = H.  $$
The vector potential $\vec A=(Hz,0)$ is uniquely determined
by the gauge given by $\BPP(N,s,q,\vec k)$, and integration
of the second relation leads to the solution space described by
 (\ref{minzero}) above, with $\alpha_n\in\RR$
undetermined constants of integration.  Therefore the
question reduces to determine when the family (\ref{minzero})
 belongs to $\BPP(N,s,q,\vec k)$.
Applying condition (\ref{phi}) we obtain
$$  2Hpqn - \omega = 2\pi k_n, \quad \mbox{for all $n\in \ZZ$.} $$
Since we have defined $k_0=0$, it follows that
$\omega=0$ and therefore
 $Hpq=m\pi$ and $k_n=m n$ for some $m\in \ZZ$.
The exact same argument applies to prove the analogous statement
in (a).

To prove (b), first observe that the upper bound on the
infimum follows from inserting the test configuration
$\psi_n$, $\vec A$ above into $\GGbp$.  Then we use the
simple inequalities,
$$   \left|{d\over dx} |\psi_n|\right|^2(x) \le 
     \left|\left({d\over dx} - i A_x(x,z_n)\right)\psi_n\right|^2(x),
\qquad (1-|\psi_n(x)|)^2 \le (1-|\psi_n(x)|^2)^2,  $$
to obtain
\beann
p\sum_{n=1}^N \| 1-|\psi_n| \|^2_{H^1[x_n,x_n+2q]}
    &\le &
 \kappa^2 p\sum_{n=1}^N \int_{x_n}^{x_n+2q}
       \left\{ (1-|\psi_n|^2)^2 
             + {1\over\kappa^2} \left|\left({d\over dx} 
                  - i A_x(\cdot,z_n)\right)\psi_n\right|^2 \right\}\, dx \\
&\le& {4\pi\over H_c^2}\GGbp (\psi_n,\vec A)\le 2Nqpr.
\eeann
By the Sobolev embedding theorem in one dimension 
we conclude that for $r$ small enough there exists a constant
$C$ with $|\psi_n(x)|\ge 1- Cr^{1/2}$ and (b) is proven.
\QED

\medskip

Lemma~\ref{inf_per}  justifies certain restrictions on the family of spaces
$\BPP(N,s,q,\vec k)$ when seeking the absolute minimizers of energy.  In
particular, in the remainder of the proof we fix
$$  q=q_m:=m{\pi\over Hp}, \ m\in {\bf N},
\qquad \vec k = m\, \ZZ.  $$
Furthermore, statement (b) justifies the treatment of the
problem in terms of the polar coordinates
$\psi_n(x)=f_n(x)\exp(i\phi_n(x))$, which in turn allows us
to restrict the winding numbers $k_n$ as above.
In the space $\BPPP$
the applied flux per period per plane is $2\pi m$, and
from (\ref{flux1}) the total flux through a period cell is
$$\iint_\Pi \curl \vec A\, dz\, dx = 2\pi Nm.  $$
Therefore, for any configuration in $\BPPP$
the  mean magnetic field in a period cell coincides with
the applied magnetic field:
$$    \meanh=
  {1\over 2q_m Np}\iint_\Pi \curl \vec A\, dz\, dx = H.  $$
This suggests that we are indeed in the ``transparent'' state
referred to in the Introduction.

Note in addition that the choice $m=1$ gives
 the smallest $q$ for which the minimum
energy is $O(r)$ for $r$ small, and $2q_1$ corresponds to the
{\em minimal period} of the explicit minimizer (\ref{minzero})
of the $r=0$ problem. Furthermore, 
the condition $\vec k=\ZZ$ will imply that
there is one Josephson vortex in between each adjacent pair
of SC planes, per period in $x$.  We will see that
the minimizing configuration in $\BPPP$ will
in fact have {\it minimal} period $2q_1$, and hence
reside in the space $\BPP(N,s,q_1,\ZZ)$.

\medskip

In summary, with the above restrictions on the period and
winding numbers, the spaces of interest are
$$   \BPP\left(N,s,{m\pi\over Hp}, m\, \ZZ\right).  $$
We recall that this space is defined as all $(f_n,\phi_n,\vec A)$ satisfying
(\ref{BP1}), (\ref{BP5}), (\ref{BP6}), and (\ref {P4}) with
$q=q_m$ and $\meanh=H$, and
\bea\label{newphix}
&&  \phi_n(x+2q_m)=\phi_n(x) + 2\pi mn;  \\
\label{newphiz}
&&  \phi_{n+N}(x+s) 
       = \phi_n(x) + NpH\, x + d.
\eea

\bigskip

\noindent
{\bf  Identifying the solution set at ${\bf r=0}$.}\
Having restricted our  choice of period and winding numbers,
thanks to Lemma~\ref{inf_per}
we may now identify the
manifold of all minimizers to the $r=0$ problem, and
verify that the setting is appropriate to apply the
degenerate perturbation theory as developed in
\cite{ACE}, \cite{AB}, \cite{ABB1}.  
\begin{prop}\label{per_r=0}  
Assume $q=q_m$, $\vec k = m\in \ZZ$ for some $m\in \NN$, and
$r=0$.
\begin{enumerate}
\item[(a)] $\inf \{\Ombpz(f_n,\phi_n,\vec A): \
   (f_n,\phi_n,\vec A)\in \BPPP\} \, =\, 0.$  The
minimum value
   is attained, and the set of all minimizers is the
$(N-1)$-dimensional hyperplane
   \bea\label{permin0}
      \SSS & := &\{(f_n,\phi_n,\vec A)\in \BPPP: \ f_n\equiv 1, \  
      \phi_n(x)=\alpha_n +nHpx, \\
       & & \qquad   \  A_x=Hz, \ A_z=0, \ 
       \mbox{where $\alpha_0=0=\alpha_1$,\ 
          $(\alpha_2,\dots,\alpha_N)\in \RR^{N-1}$.}  \} \nnn
     \eea
\item[(b)]  For any element $\sigma=(f_n^0,\phi_n^0,\vec A^0)\in \SSS$, the
     linearized operator
    $D^2\Ombpz(\sigma): \ E\to E$ defines a Fredholm operator with index
zero.
      Moreover,
 \be\label{pernullspace}
    T_\sigma\SSS=\ker D^2\Ombpz(\sigma)\simeq \RR^{N-1}. \ee
\end{enumerate}
\end{prop}
Recall that the space 
$E=E(N,s,{m\pi\over pH})$ is the tangent space to
$\BPPP$, and was defined at the end of section~3.
We define the second variation of energy as a quadratic
form on $E$:  in particular, for $(f_n^0,\phi_n^0,\vec A^0)\in \SSS$
and $(u_n,v_n,\vec a)\in E$, 
\bea \nnn
D^2\Om_0 (f_n^0),\phi_n^0,\vec A^0)[u_n,v_n,\vec a] &=&  
  p\sum_{n=1}^N \int_{x_n-q_m}^{x_n+q_m} \left\{ 2u_n^2 +
      {1\over\kappa^2}  [u'_n]^2 + {1\over\kappa^2}[v'_n-a_x(x,z_n)]^2
        \right\}\, dx  \\
   \nnn
   && \qquad
   +  {1\over\kappa^2}\iint_{\Pi}\left| \curl \vec a \right|^2 \, dx\, dz .
\eea
The proof of Proposition~\ref{per_r=0} (a)
follows easily from the proof of Lemma~\ref{inf_per}~(c), and
part (b) is clear from the form of $D^2\Om_0 (f_n^0),\phi_n^0,\vec A^0)$
above (see Proposition~2.1 in \cite{ABB1} for details.)

\medskip

Note that $\TT= T_\sigma\SSS$ is independent of
$\sigma\in\SSS$, and that we may
treat $\SSS$ as a (compact) $(N-1)$-torus
${\bf T}^{N-1}$, since the  energy is $2\pi$-periodic in each $\alpha_n$.
 Applying the  condition (\ref{BP4}) 
to an element of $\SSS$  we obtain the periodicity
conditions for $\alpha_n$ and $\delta_n:=\alpha_n-\alpha_{n-1}$:
\be\label{alpha_n}
\alpha_{N+n} -\alpha_n =d-Hps(n+N), \quad
   \delta_{n+N}-\delta_n = -Hps \ \mbox{(mod $2\pi$)},
\ee
where $d$ is the constant which appears
in the 't Hooft-type periodicity condition (\ref{newphiz}).
 Since we are given $\alpha_0=0=\alpha_1$ in the definition
 of the space $\BPPP$,
(\ref{alpha_n}) provides the relation $\alpha_N=d-HNps$. 
This reveals how a choice of $(\alpha_2,\dots,\alpha_N)$
determines all $\alpha_n$:  $\alpha_N$ determines the
constant $d$, and the first equation of (\ref{alpha_n}) 
then generates all the others in the sequence. 

\smallskip

We also note that the $\alpha_n$ ($n=2,\dots,N$)
may be used as parameters for the manifold $\SSS$.
However, it will be more convenient in the end to
parametrize $\SSS$ by the phase differences,
$$   \delta_n:= \alpha_n-\alpha_{n-1}.  $$
Since the values of $(\delta_2,\dots,\delta_N)$  determine
$(\alpha_2,\dots,\alpha_N)$ uniquely, this is an equivalent
choice of parametrization.
We abuse notation and denote 
$$   \sigma=\sigma(\delta_2,\dots,\delta_N)\in \SSS,  $$
$(\delta_2,\dots,\delta_N)\in\RR^{N-1}$.

\section{Degenerate perturbation theory.}
\setcounter{thm}{0}

We now perturb away from the degenerate minima of
$\Ombpz$, using a variational
Lyapunov--Schmidt procedure,
just as in \cite{ABB1}.  This method
has been used by Ambrosetti, Coti-Zelati, \& Ekeland \cite{ACE},
 Abrosetti \& Badiale \cite{AB}, Li \& Nirenberg \cite{LN} (and
 many others) in a variety of
situations
 involving heteroclinic solutions of Hamiltonian systems and
 in the semiclassical limit of the nonlinear Schr\"odinger equation.

\medskip

Since $\SSS$ is a hyperplane, $\TT=T_\sigma\SSS$ is independent of
$\sigma\in\SSS$.  Let $W=\TT^\perp$, so any 
$(f_n,\phi_n,\vec A)\in \BPPP$
admits the unique decomposition 
$(f_n,\phi_n,\vec A)=\sigma+w$ with $\sigma\in\SSS$,
$w\in W$, and any $U:=(u_n,v_n,\vec a)\in E$ decomposes uniquely
as $U=t+w$ with $t\in\TT$, $w\in W$.
We denote the orthogonal projection maps $P: E\to \TT$, $P^\perp:E\to W$
so that $PU=t$, $P^\perp U=w$ when $U=t+w$.   
We exploit the Hilbert space setting and interpret the first variation
$\nabla\Ombp(f_n,\phi_n,\vec A)$ as an element of $E$ itself, and
 project the equation $\nabla\Ombp(f_n,\phi_n,\vec A)=0$
into the two linear subspaces $\TT$ and $W$,
\bea
\label{Teqn}
F_1(r,\sigma,w) & :=  & P\left[ \nabla\Ombp (\sigma+w) \right] = 0; \\
\label{Weqn}
F_2(r,\sigma,w) & :=  & P^\perp\left[ \nabla\Ombp (\sigma+w) \right] = 0.
\eea
The second equation can be solved uniquely for
$w=w(r,\sigma)$ in a neighborhood
of $\SSS$ for $r$ small, using the Implicit Function Theorem.
Because our functional $\Ombp$ is smooth we can expand 
$w(r,\sigma)$ in powers of $r$, and since
$\Ombp(\sigma+w)$ is periodic in 
$\sigma$ the expansion is uniform $\sigma$.   The
result below is based on Lemma~2 of \cite{AB}, and
follows directly from Lemma 3.1, \cite{ABB1}:
\begin{lem}\label{Lya}  Assume $q=q_m={m\pi\over Hp}$,
$\vec k=m\ZZ$ for some $m\in \NN$.

There exist constants $r_0>0$ and $\delta>0$,
depending on $N,m,\kappa,H,s$  and a
smooth function
$$  w=w(r,\sigma): \ (-r_0,r_0)\times \SSS\to W\subset E  $$
such that:
\begin{enumerate}
\item[(i)] There exists smooth functions
$w_1,w_2$
such that 
$$  w(r,\sigma)=r w_1(\sigma) + r^2 w_2(r,\sigma)
$$
 for all $|r|<r_0$ and for all $\sigma\in\SSS$;
\item[(ii)]   $P^\perp [\nabla \Ombp(\sigma+w(r,\sigma))]=0$.  
\item[(iii)]
Conversely,
if $P^\perp [\nabla \Ombp(\sigma+w)]=0$ for some
$r\in(-r_0,r_0)$ and  $w\in W$ with $\|w\|_E < \delta$, then
$w=w(r,\sigma)$.
\item[(iv)]
For any choice of $m_0,\kappa_0,s_0, H_0>0$ the constant $r_0$
may be chosen uniformly for all $N\ge 1$, $1\le m \le m_0$,
$1\le \kappa\le \kappa_0$, $H\ge H_0$ and $|s|\le s_0$.
\end{enumerate}
\end{lem}
Parts (i)--(iii) follow easily from the Implicit
Function Theorem.  The dependences on the various parameters  in (iv)
is more delicate, and was the subject of section~5 of \cite{ABB1}.
In that section we proved a lower bound on $r_0$ via {\it a priori}
estimates on the solutions of (\ref{Weqn}). 
We note that  the 
interval of validity of the expansion in \cite{ABB1} was strongly affected
by the sample width $L$, with improved convergence with smaller $L$.
In the periodic problem the same method as in 
section~5 of \cite{ABB1} can be used to obtain the same lower bound on
the radius of validity, but with the period $q=q_m$ replacing $L$ in
the expression of the lower bound.  Since $2q_m$ {\it decreases}
with increasing field strength $H$ we can expect our solutions
to have a large range of validity in high fields.

\medskip

We define 
$$  \SSS_r := \{\sigma +w(r,\sigma): \ \sigma\in \SSS\}.  $$
$\SSS_r$ is a smooth manifold parametrized by the
hyperplane $\SSS$.  The important role played by
$\SSS_r$ is that it is a natural constraint for $\Ombp$
(see Lemma~4 of \cite{AB}),
and hence the equation (\ref{Teqn}) may be solved variationally:
\begin{lem}\label{natural}  Assume $q=q_m={m\pi\over Hp}$,
$\vec k=m\ZZ$ for some $m\in \NN$.
\begin{enumerate}
\item[(a)]\
If $(f_n,\phi_n,\vec A)\in \SSS_r$ satisfies
$D(\Ombp{}_{|\SSS_r})(f_n,\phi_n,\vec A)=0$, then
$\nabla\Ombp(f_n,\phi_n,\vec A)=0$ in $E$.
\item[(b)]\
There exists a constant
$\tilde r_0=\tilde r_0(N,m,\kappa,s,H)$ with
$0<\tilde r_0< r_0$ such that if
$(f_n,\phi_n,\vec A)\in \BPPP$ is a critical point
of $\Ombp$ with 
$$ \Ombp(f_n,\phi_n,\vec A)\le 2q_mNp\, r = {2\pi  mN\over H}\, r
       \quad\mbox{and}\quad
          |r|<\tilde r_0,   $$
 then $(f_n,\phi_n,\vec A)\in\SSS_r$.
\end{enumerate}
\end{lem}
The consequence of this lemma is very important:  
for $r$ sufficiently small the absolute minimizer will
be found by minimization on
 the finite dimensional manifold $\SSS_r$.
We note that the parameter dependences for the interval
of validity $r_0$ in Lemma~\ref{Lya} do not carry
through to Lemma~\ref{natural}.  This is because we
have no estimate on the neighborhood $\delta$ of
Lemma~\ref{Lya} in terms of other parameters,
and so we cannot  conclude that configurations with very small
energy must be ``close enough'' to $\SSS$ to be in the range
of Lemma~\ref{Lya}~(iii).  

 The proof of Lemma~\ref{natural} is identical to 
Lemma~3.2 in \cite{ABB1} and is omitted.

\section{Vortex lattice solutions.}
\setcounter{thm}{0}

We are now ready to treat the finite dimensional problem
(\ref{Teqn}), via minimization of the constrained functional
$\Ombp|_{\SSS_r}$.  As opposed to the finite-width case
studied in \cite{ABB1} the functional $\Ombp|_{\SSS_r}$
will degenerate at order $r$, and therefore it will be
necessary to carry out the expansion to higher order.  
To accomplish this we must calculate an expansion of
the solutions of the regular projected equation (\ref{Weqn})
to order $r$.  Fortunately, the regularity of $\Ombp$ and
the Implicit Function Theorem allow us to do the computation
explicitly, and we will obtain a straight-forward expansion
of $\Ombp|_{\SSS_r}$ which permits direct minimization.

We summarize our conclusion in the following Theorem:
\begin{thm}\label{Thm3}  Let $q_m={m\pi\over Hp}$, $\vec k=m\ZZ$,
$m\in\NN$.

\noindent
\ (i)\
For every  $s\in \RR$ and $m=1,2,\dots$,
there exists $r_0=r_0(N,s,\kappa,m,H)>0$ 
such that for all $r\in (0,r_0)$ the minimizer of $\Ombp$ in $\BPPP$
is a $\Pi_{N,s,q_m}$-periodic solution
given asymptotically  by 
\be \label{XXXX}  \left.  \begin{array}{c}
f_n =  1\, +\, r\left[ 
-\frac12 + {\kappa^2\over 2(H^2p^2 + 2\kappa^2)}
      \left( \cos(\delta_n + Hp x) + \cos(\delta_{n+1} + Hp x)\right)\right]
              +o(r),  \\
\\
 h(x,z)=H -r{\kappa^2\over 2H}\cos(\delta_n+Hpx) + o(r)\\
\\
j_x^{(n)}(x)=r{\kappa^2\over 2Hp}\left(\cos(\delta_{n+1}+Hpx) 
-\cos(\delta_{n}+Hpx)\right)+ o(r)\\
\\
j_z^{(n)}(x)= {1\over 2}\, r\, \kappa^2 p\sin (\delta_{n}+Hpx) + o(r),
\end{array} \right\}
\ee
 where
$(\delta_2,\dots,\delta_N)$ is a minimizer of the
finite dimensional problem
\be\label{FDprob} 
  F(N,s):= 
 \inf
\left\{ {1\over N}\sum_{n=1}^N \cos(\delta_n - \delta_{n+1}): \ 
(\delta_2,\dots,\delta_N)\in \RR^{N-1}, \ \delta_1=0,
\ \delta_{N+1}=-Hps
\right\}.
\ee
Moreover, the minimum energy is given by:
\be\label{minenergy}  \inf_{\BPPP}\Ombp =
    2mq_1Np\left( r + r^2 \left(  C_0 + C_1 F(N,s) \right)\right)
   + O(r^3),  
\ee
where $q_1={\pi\over Hp}$ and 
$C_0\in\RR$, $C_1>0$ are constants independent of $N,s,q,m$.

\noindent
(ii)\quad
 $\inf\{ F(N,s): \ s\in\RR, \ N=1,2,\dots\}=-1$, \quad
 and the minimum is attained at $(N,s)$ 
if and only if  either:
$$ \mbox{$N$ is even and $s=2\ell q_1$ for $\ell\in \ZZ$;}  $$
or if 
$$\mbox{$N$ is odd and $s=(2\ell+1) q_1$ for $\ell\in \ZZ$.}  $$
In either case, the minimizer of $\Ombp$ in $\BPPP$ is {\em unique} and
coincides with the period-$2p$ in $z$,
period $2q_1$ in $x$
lattice which minimizes $\Ombp$ in $\BPP(1,q_1,q_1,\ZZ)$,
satisfying (\ref{shiftper}) and given asymptotically by
(\ref{per_2}).
\end{thm}

We observe that the independence of the constants $C_0,C_1$
with respect to $N,s,q,m$ means that the dependence of the geometry
of $\Pi$ on the energy per cross-sectional
unit area is entirely encoded in the function $F(N,s)$.
In the generic case (i), $F(N,s)$ could have several absolute minimizers
$(\delta_2,\dots,\delta_N)$,
and choosing different sequences of $r\to 0$ could lead to
minimizers of $\Ombp|_{\SSS_r}$ which are determined by
different minimizers of $F(N,s)$.  In the special case (ii)
the absolute minimizer of $F(N,s)$ is unique (and non-degenerate)
and we obtain the smallest possible energy for the Lawrence--Doniach
energy, corresponding to the second
alternative in Theorem~\ref{bigthm}

\medskip

\noindent
{\it  Proof of Theorem~6.1:}\
By Lemma~\ref{Lya}, we may decompose an element of $\BPPP$ as
$(f_n,\phi_n,\vec A)= \sigma+ w(r,\sigma)$, with $\sigma\in\SSS$ and
$w\in W=[T\SSS]^\perp\subset E$ such that the Euler-Lagrange equations
hold when projected into the space $W$.   Furthermore, we may
write $w(r,\sigma)= rw_1(\sigma) + r^2 w_2(\sigma,r)$,
and in $(f_n,\phi_n, A_x, A_z)$ coordinates, we denote
$$   w_1=(u_{n,1},v_{n,1},a_{x,1},a_{z,1}).
$$
In other words, recalling (\ref{permin0}),
\be \label{wexpansion}  \left. \begin{array}{lll}
   f_n=1+ru_{n,1} + O(r^2),  &\quad  & 
     \phi_n= \alpha_n + nHpx + v_{n,1}+O(r^2),  \\
           A_x= Hz + ra_{x,1} + O(r^2), &\quad  &
            A_z=ra_{z,1}+O(r^2).
   \end{array}  \right\}
\ee
Note also that $w_1(\sigma)=\partial_r w(0,\sigma)$.

\smallskip
\noindent
{\bf Step 1:}\ Expansion of the energy.

We recognize that the energy can be written in two
parts, $\Ombp(U)=\Ombpz(U) + r G(U)$.
Since $\Ombp|_{\SSS_r}=\Ombp(\sigma + w(r,\sigma))$ is a smooth
function of $r$ and $\sigma\in\SSS$,
it admits a Taylor expansion at $r=0$ of the form,
$$
\Ombp(\sigma+w(r,\sigma))= \Omega^{(0)} + r\Omega^{(1)}
     + {r^2\over 2} \Omega^{(2)} + O(r^3) ,
$$
with $O(r^3)$ remainder uniform in $\sigma\in\SSS$, where:
\begin{eqnarray*}
\Omega^{(0)} &=& \Ombpz(\sigma)=0, \\
\Omega^{(1)} & =& \left.  {d\over dr}\right|_{r=0}
                    \Ombp(\sigma + w(r,\sigma)) \\
&=&  G(\sigma) + \nabla\Ombpz(\sigma)[w_1(\sigma)] = G(\sigma) \\
 &=&  p \sum_{n=1}^N \int_{x_n}^{x_n+2q_m}
        (1-\cos(\delta_n + Hp x))\, dx  \\
&=& 2Npq_m, \\
\Omega^{(2)} & = &  \left. {d^2\over dr^2}\right|_{r=0}
                    \Ombp(\sigma + w(r,\sigma)) \\
&=& 2\nabla G(\sigma)[w_1(\sigma)] +
        \nabla^2 \Ombpz (\sigma) [w_1,w_1]
           + \nabla\Ombpz(\sigma)\left[{\partial^2 w\over \partial r^2}\right] \\
& =&  2\nabla G(\sigma)[w_1] + \nabla^2\Ombpz(\sigma)[w_1,w_1] \\
&=& p \sum_{n=1}^N \int_{x_n}^{x_n+2q_m}
     \left\{   2u_{n,1}^2 + {1\over \kappa^2} (u'_{n,1})^2
          +  {1\over \kappa^2} (v'_{n,1} - a_{x,1}( x,z_n))^2 
          \right\}\, dx  \\
    & &   + 
   p \sum_{n=1}^N \int_{x_n}^{x_n+2q_m}
     \left\{  (u_{n,1} + u_{n-1,1})[1-\cos(\delta_n + Hp x)]
         + \sin(\delta_n +Hp x) \varphi_{n,n-1,1}(x)  \right\}\, dx  \\
         && 
     + {1\over \kappa^2}
       \iint_\Pi
                \left[ {\partial a_{x,1}\over\partial z}  -
                    {\partial a_{z,1}\over\partial x} \right]^2 \, dz \, dx.
\end{eqnarray*}
Here we denote 
\be\label{varphi}
 \varphi_{n,n-1,1}=v_{n,1} - v_{n-1,1} 
   - \int_{z_{n-1}}^{ z_n} a_{z,1}(x,z)\, dz.
\ee
Note that the term of order $r$ degenerates-- unlike the
boundary-value problem treated in \cite{ABB1} the 
constant phase differences $\delta_n$ are
{\it not} determined at this point, but only at the order $r^2$!

\smallskip
\noindent
{\bf Step 2:}\  Expansion of the solutions.

To evaluate the next order term $\Omega^{(2)}$ we require
the first-order correction to the solution
 $w_1(\sigma)$.  Implicit
differentiation of the equation (\ref{Weqn}) with respect to $r$ yields:
\be\label{firstorder}
 P^\perp\left[ \nabla G(\sigma) + \nabla^2\Ombpz (\sigma)w_1
    \right] =0.  
\ee
Since $\nabla^2\Ombpz(\sigma)$ is invertible on $W$ (\ref{firstorder})
determines $w_1=w_1(\sigma)$ uniquely.

Using the expansion (\ref{wexpansion})
we write the projected equation (\ref{firstorder}) in terms of the above
coordinates:
\bea
\label{proj1}  && 
-{1\over \kappa^2} u''_{n,1} + 2u_{n,1}  = 
          {1\over 2}\left(
       \cos\Phi_{n,n-1,0}
      +\cos\Phi_{n+1,n,0}-2\right), \\
\label{proj2} &&  
{1\over \kappa^2}{d\over dx}\left( v'_{n,1} - a_{x,1}(x,z_n)\right) = 
       {1\over 2}[\sin(\delta_n + Hpx)
       - \sin(\delta_{n+1}+Hpx)], 
\eea
and
\be
\label{proj3} 
\vec a_1(x,z)= \left({\partial\xi\over\partial z}, 
       -{\partial\xi\over\partial x}\right),
          \qquad \Delta\xi=b_1(x,z), 
              \quad \left.\xi \right|_{\partial B}=0,
\ee
where
\be
\label{proj4} \left.
\begin{array}{c}
   b_1(x,z)=b^{(n)}_1(x), \quad z_{n-1}<z<z_n, n=1,\dots,N, \\
   \\
   b^{(n)}_1(x)-b^{(n+1)}_1(x) = p(v'_{n,1} - a_{x,1}(x,z_n)), n=1,\dots,N-1, \\
   \\
  \partial_x b^{(n)}_1=
          \frac12\, p\kappa^2\, \sin(\delta_n + Hpx), \\
        \\
   b^{(n)}_1(\pm L)=0.
\end{array} \right\}
\ee
 This system also coincides with
the Euler--Lagrange equations for minimizing
$\Omega^{(2)}$ above in the space $W$.
An equivalent way to arrive at these
equations is to begin with the equations for $w(r,\sigma)$
satisfying (\ref{Weqn}) as derived in
 section~2.5 of \cite{ABB1}, then take the order $r$ terms
appearing in each equation.

\medskip

These equations may be integrated explicitly:
from the current conservation equation 
(\ref{proj2}) we obtain:
\be\label{bp1}
v_{n,1}' - a_{x,1}(x,z_n) = 
C_n - {\kappa^2\over 2Hp} 
   \cos(\delta_n + Hpx)
      + {\kappa^2\over 2Hp} \cos(\delta_{n+1} + Hpx), \quad n=1,\dots N.
\ee
where $C_n$ are as-yet undetermined constants.
Next we use the equations (\ref{proj4}) for the magnetic
field inside each gap, $b_1^{(n)}(x)=b_1(x,z)$ for $z_{n-1}<z<z_n$,
\be\label{p2}
b_1^{(n)}(x)= -{\kappa^2 \over 2H} \cos(\delta_n+Hpx) + D_n,
    \quad n=1,\dots,N,
\ee
with $D_n$ another set of undetermined constants.
Using the jump condition of (\ref{proj4}) for the magnetic
field across each superconducting plane and the periodicity 
of $D_n$, we see that
\bea\label{bp3}
-pC_n = (D_{n+1}-D_n),\qquad n=1,\dots,N.
\eea
We now solve for the order $r$ term in the
gauge-invariant phase $\varphi_{n,n-1,1}$ (defined in (\ref{varphi}))
 using formula (\ref{Phieqn}).
We find for $n=1,\dots, N$,
\beann
{d\over dx}\varphi_{n,n-1,1}(x)
  &=&  (v'_{n,1}-a_{x,1}(x,z_n)) -(v'_{n-1,1}-a_{x,1}(x,z_{n-1})) +p  b_1^{(n)}(x)
\\
   &=& pD_n + C_n - C_{n-1} \\
   && +{\kappa^2\over 2Hp} \left( \cos(\delta_{n+1} + Hpx)
          -(2+p^2)\cos(\delta_n + Hpx) + \cos(\delta_{n-1} + Hpx) \right).
\eeann  
Since each $\varphi_{n,n-1,1}(x)$ must be $2q$-periodic in $x$ we
must satisfy the integrability conditions $pD_n + C_n - C_{n-1}=0$.
Together with (\ref{bp3}), we
derive a second order difference equation for $D_n$,
\be
\label{difference_eqn}
   D_{n+1} -2 D_n + D_{n-1} = p^2 D_n, \qquad n=1,\dots,N,  
\ee
and note that the $\vec e_2$-periodicity condition (\ref{hper}) of $h$ 
together with the corresponding condition (\ref{alpha_n}) for
$\delta_n$
implies $D_n=D_{n+N}$.
The maximum principle for second order difference equations
then ensures that $D_n=0$, and hence  $C_n=0$ by (\ref{bp3}),
and all constants are uniquely determined.  We may then integrate
to obtain:
$$  \varphi_{n,n-1,1}(x)
 ={\kappa^2\over 2H^2 p^2}  \left( \sin(\delta_{n+1} + Hpx)
          -(2+p^2)\sin(\delta_n + Hpx) + \sin(\delta_{n-1} + Hpx) \right).
 $$
We note that the arbitrary constant of integration which should
normally come with $\varphi_{n,n-1,1}$ is zero here, since we are
solving for each  $\phi_n$ in the orthogonal to $T\SSS$.
(That constant would result from the order-$r$ correction to our
choice of $\sigma\in \SSS$ in the minimization problem on
$\SSS_r$.)

 In conclusion, the gauge-invariant quantities associated with
points $(\sigma+w(r,\sigma))$ on the natural constraint $\SSS_r$ are:
\be \label{persol}  \left.  \begin{array}{c}
f_n =  1\, +\, r\left[ 
-\frac12 + {\kappa^2\over 2(H^2p^2 + 2\kappa^2)}
      \left( \cos(\delta_n + Hp x) + \cos(\delta_{n+1} + Hp x)\right)\right]
              +O(r^2),  \\
\\
 h(x,z)=H -r{\kappa^2\over 2H}\cos(\delta_n+Hpx) + O(r^2)\\
\\
j_x^{(n)}(x)=r{\kappa^2\over 2Hp}\left(\cos(\delta_{n+1}+Hpx) 
-\cos(\delta_{n}+Hpx)\right)+ O(r^2)\\
\\
j_z^{(n)}(x)= {r\over 2}\kappa^2 p\sin (\delta_{n}+Hpx) + O(r^2),
\end{array} \right\}
\ee
where (as usual) $\sigma=\sigma(\delta_2,\dots,\delta_N)$.
We emphasize that  Lemma~\ref{Lya}~(i)
ensures that all remainder terms are uniform in $\sigma\in \SSS$.

\medskip

\noindent
{\bf Step 3:}\   Expansion of $\Ombp|_{\SSS_r}$.

We now resolve the degeneracy at order $r^2$ to determine which
choice of $\sigma=\sigma(\delta_2,\dots,\delta_N)$ in $\SSS$ gives rise
to stationary solutions of the Lawrence--Doniach system.
Substituting and computing the integrals,
\beann
\Omega^{(2)} & =& 
 \left\{  Npq_m {\kappa^2\over 2H^2p^2}
   \left[ \frac32 p^2 -1 - {H^2p^2\over H^2p^2 + 2\kappa^2} \right]
        - Npq_m   \right\}  \\
     & &  \quad
     +\  q_m {\kappa^2\over 2H^2p^2}
         \left( 1-{H^2p^2\over H^2p^2 +
            2\kappa^2}\right) \,
             p\sum_{n=1}^N  \cos(\delta_n - \delta_{n+1}) ,
\eeann
where we recall $q_m=m\pi/Hp$.

In conclusion, we obtain the following expansion
of $\Ombp|_{\SSS_r}$, with $\SSS_r$ parametrized
by $\sigma=\sigma(\delta_2,\dots,\delta_N)$:
\be\label{energy2} \Ombp(\sigma + w(r,\sigma)) = 2Npq_m\left\{ r 
   + r^2\left( C_0 + C_1\, {1\over N}\sum_{n=1}^N  
      \cos(\delta_n - \delta_{n+1}) 
           \right) \right\}
       + O(r^3) ,
\ee
where $C_0\in\RR$, $C_1>0$ are  constants 
independent of $N,s,q,m$.  The periodicity
conditions in $\BPPP$ carry over to an inhomogeneous
boundary condition on $\delta_n$ (see
(\ref{alpha_n})),
\be\label{bc} \delta_1=0, \qquad \delta_{N+1} = -Hps \ \mbox{(mod $2\pi$).}  
\ee
We recall once again that the remainder term is {\it uniform}
in $\sigma=\sigma(\delta_2,\dots,\delta_N)\in \SSS$.

\bigskip

\noindent
{\bf Step 4:}\  Verifying equation (\ref{minenergy}), and the
conclusion of (i).

First we observe that by Lemma~\ref{inf_per} and Lemma~\ref{natural}~(b)
for $0<r<\tilde r_0$ the infimum of $\Ombp$ in $\BPPP$ will
be attained on $\SSS_r$.
Since $\SSS_r$ is finite dimensional
and $\Ombp|_{\SSS_r}$ is periodic in the local coordinates
$(\delta_2,\dots,\delta_N)$,  for every
$0<r<\tilde r_0$ there exists (at least one) minimizer, 
$$   \Ombp(\sigma_r + w(r,\sigma_r))=\inf \Ombp|_{\SSS_r} , \qquad 
       \sigma_r=\sigma(\delta_2(r),\dots,\delta_N(r))\in \SSS.
$$
By the expansion (\ref{energy2}) we have
\bea\nnn
 \Ombp(\sigma_r+ w(r,\sigma_r)) &=&
2Npq_m\left\{ r 
   + r^2\left( C_0 + C_1\, {1\over N}\sum_{n=1}^N  
      \cos(\delta_n(r) - \delta_{n+1}(r)) 
           \right) \right\}
       + O(r^3) \\
\label{lowerbound}
&\ge & 2Npq_m\left\{ r 
   + r^2\left( C_0 + C_1\, F(N,s) 
           \right) \right\}
       + O(r^3),
\eea
since $F(N,s)$ is the infimum of the sum of cosines over
all possible configurations.  To obtain a complementary
inequality, let $(\delta^*_2,\dots,\delta^*_N)$ be
any minimizer of $F(N,s)$, that is
$$   F(N,s)={1\over N}\sum_{n=1}^N \cos (\delta^*_n - \delta^*_{n+1}).  $$
(Under the hypotheses of (i) there could be many such minimizers.)
Then, applying (\ref{energy2}) to this configuration, we obtain:
\bea\nnn   \inf \Ombp|_{\SSS_r}&\le &
        \Ombp(\sigma^* + w(r,\sigma^*))
\\
\nnn
&\le& 2Npq_m\left\{ r 
   + r^2\left( C_0 + C_1\, {1\over N}\sum_{n=1}^N  
      \cos(\delta^*_n - \delta^*_{n+1}) 
           \right) \right\}
       + O(r^3)   \\
\label{upperbound}
& = & 2Npq_m\left\{ r 
   + r^2\left( C_0 + C_1\, F(N,s)  \right) \right\}
       + O(r^3).
\eea
Putting together (\ref{lowerbound}) and (\ref{upperbound}) we
deduce the  energy expansion (\ref{minenergy}) stated in 
Theorem~\ref{Thm3}~(i).  Moreover, (\ref{upperbound})
implies
$$  F(N,s)\le {1\over N}\sum_{n=1}^N  
      \cos(\delta_n(r) - \delta_{n+1}(r)) 
             \le F(N,s) + O(r).  $$
Therefore for any sequence of $r\to 0$, 
 the minimizers $\sigma_r$ of $\Ombp|_{\SSS_r}$ form a 
minimizing sequence for the variational problem $F(N,s)$.
Hence the $\sigma_r$
accumulate as $r\to 0$ at minimizers of $F(N,s)$.  To be more
precise, for any sequence of $r\to 0$ there exist subsequences
and minimizers $\sigma^*=\sigma(\delta^*_2,\dots,\delta_N^*)$ of $F(N,s)$
such that (along the subsequence) $\sigma_r\to\sigma^*$.
Inserting this information into (\ref{persol}) we
obtain (\ref{XXXX}).  This completes the proof of part~(i)
of Theorem~\ref{Thm3}.

\medskip

\noindent
{\bf Step 5:}\ Proof of (ii).

  By the
expansion (\ref{minenergy}) the problem reduces
to determining for which lattice parameters $N,s$ 
does $F(N,s)$ attain its
lower bound of $-1$.  This lower bound is achieved
if and only if the boundary condition
(\ref{bc}) admits a choice of $\delta_n$ with $\delta_{n+1}-\delta_n=\pi$
(mod $2\pi$.)  

\smallskip

When $s\not\in {\pi\over Hp}\ZZ$ the lattice is
frustrated since the space $\BPPP$ does not admit the configuration
with $\delta_{n+1}-\delta_n=\pi$ (mod $2\pi$.)  In that case
we obtain $F(N,s)>-1$, and the energy per unit area of the
minimizer in $\BPPP$ will be strictly larger than
this absolute minimum value, for all sufficiently small
$r>0$.

\smallskip

If 
\be\label{even}
\mbox{$N$ is even and $s=2\ell {\pi\over Hp}$ for $\ell\in \ZZ$}
\ee
or if 
\be\label{odd}
\mbox{$N$ is odd and $s=(2\ell+1) {\pi\over Hp}$ for $\ell\in \ZZ$,}
\ee
the choice $\delta_n=(n-1)\pi$ (mod $2\pi$) is allowed by 
(\ref{bc}) and the infimum $F(N,s)=-1$ 
is attained.   Define $g:\RR\times \RR^{N-1}\to \RR$ with
\beann  g(r,\delta_2\dots,\delta_N)&:=& {1\over r^2} \left(
     {\Ombp(
         \sigma+ w(r,\sigma))
        \over 2mq_1 Np}-  r  \right) \\
   &=& C_0 + C_1 \frac1N \sum_{n=1}^N \cos(\delta_n - \delta_{n+1}) + O(r),
\eeann
by (\ref{energy2}).  In particular, $g$ is smooth and 
$g(0,\delta_2\dots,\delta_N)$ is minimized if and only if
$\delta_n=\delta_n^*=(n-1)\pi$~mod~$2\pi$.
Now, $(\delta_2^*,\dots,\delta^*_N)$ is a {\it non-degenerate}
minimizer of $g(0,\delta_2\dots,\delta_N)$:  its Hessian is
the familiar tridiagonal matrix with $2$ on the diagonal and
$-1$ on each off-diagonal, associated with the 
(positive-definite) second-order difference
operator.
By the Implicit Function Theorem we  conclude that 
for all sufficiently small $r>0$ the function $g(r,\delta_2\dots,\delta_N)$
has a unique minimum at $(\delta_2(r),\dots,\delta_N(r))$,
with $\delta_n(r)=\delta_n^* + O(r).$   In other words,
setting $\sigma_r:=\sigma(\delta_2(r)\dots,\delta_N(r))$,
$$   \inf \Ombp|_{\SSS_r} = \Ombp(\sigma_r + w(r,\sigma_r)).  $$
By Lemma~\ref{natural} when $0<r<\tilde r_0$ this gives
the global minimizer of $\Ombp$ in $\BPPP$, and inserting these
optimal values for $\delta_n$ into
the asymptotic formulae
(\ref{persol}) we obtain (\ref{per_2}).
This proves that the minimizers for any period geometry
satisfying (\ref{even}) or (\ref{odd}) coincide up to
order $r$.  We will do better, and show that they are
actually identical.

\smallskip

To gain a complete understanding of
the absolute lowest energy solution
 we first 
consider the special case $N=1$, $s=q_1=\pi/Hp$, and
$m=1$.  In this case, there is no degenerate manifold:  
the $r=0$ problem has a unique, non-degenerate solution and the
perturbation is regular.  In particular the second-order term
in the energy expansion (\ref{minenergy}) is completely
determined by  (\ref{bc}):   since $\delta_0=0$ and (for $N=1$)
$\delta_{N+1}=\delta_2=-Hsp=-\pi$  we have
$F(1,q_1)=-1$.  We obtain for all sufficiently small $r>0$ a
unique solution $(f^1_n,\phi^1_n,\vec A^1)$ which minimizes
$\Ombp$ in $\BPP(1,s,q_1,\ZZ)$.  The gauge-invariant quantities associated
to $(f^1_n,\phi^1_n,\vec A^1)$ will be $2q_1$-periodic in $x$
and shifting $z$ by $p$ results in a  translation
by a half-period $q_1$ in $x$ (see (\ref{shiftper}).)

\smallskip

Now return to the cases (\ref{even}) or (\ref{odd}) where
$N\ge 2$ and $F(N,s)=-1$.  
We denote by $\delta_n^*=(n-1)\pi$~mod~$2\pi$, the unique
absolute minimizer of $F(N,s)$, and by
$(\delta_2(r),\dots,\delta_N(r))$ the coordinates of the
 absolute minimizer $(\tilde f_n,\tilde\phi_n,\tilde A)\in \BPP(N,s,q_m,m\ZZ)$
of $\Ombp$, which we know
lies on ${\SSS_r}$ for all small $r$, and for which
$\delta_n(r)=\delta_n^* + O(r).$

It is
easy to verify that the $\BPP(1,q_1,q_1,\ZZ)$--minimizer
 $(f^1_n,\phi^1_n,\vec A^1)$ is also $\Pi_{N,s,q_m}$-periodic when
$N,s$ satisfy (\ref{even}) or (\ref{odd}), and in fact
 $(f^1_n,\phi^1_n,\vec A^1)$ solves the Euler--Lagrange
equations for $\Ombp$ in $\BPP(N,s,q_m, m\ZZ)$.
Moreover, $(f^1_n,\phi^1_n,\vec A^1)$ also satisfies
the expansion of energy given by (\ref{minenergy})
on the period parallelogram $\Pi_{N,s,q_m}$, and therefore
by Lemma~\ref{natural} it lies on the constraint manifold
$\SSS_r\subset\BPP(N,s,q_m,m\ZZ)$ in cases (\ref{even})
and (\ref{odd}).  Finally, the coordinates of $(f^1_n,\phi^1_n,\vec A^1)$
on $\SSS_r\subset\BPP(N,s,q_m,m\ZZ)$, which
we denote by $(\delta_1^{(1)}(r),\dots,\delta_N^{(1)}(r))$,
satisfy
$$   \delta_n^{(1)} = \delta_n^* + O(1).  $$
Since $(\delta_2^*,\dots,\delta_N^*)$ is a non-degenerate
minimizer of $F(N,s)$ we must have  $\delta_n(r)=\delta_n^{(1)}(r)$
for all sufficiently small $r$, that is 
$(\tilde f_n,\tilde\phi_n,\tilde A)=(f_n^1, \phi_n^1, \vec A^1)$
for all small $r$, that is the minimizers in $\BPP(N,s,q_m,m\ZZ)$
with (\ref{even}) or (\ref{odd}) coincide exactly with the
minimizers in $\BPP(1,{\pi\over Hp},{\pi\over Hp},\ZZ)$.

This concludes the proof of
Theorem~\ref{Thm3}. 
\QED

Finally we prove Theorem~\ref{bigthm}.
In case $q\not\in {\pi\over Hp}\ZZ$ or $k_n\neq mn$ 
(for $m\in {\bf N}$ constant) we have by Lemma~\ref{inf_per}
$$ \epsilon_r(N,s,q,\vec k)  :=  \inf_{\BPP(N,s,q,\vec k)}\Ombp 
  \ge
      \inf_{\BPP(N,s,q,\vec k)}\Ombpz =: \omega_0(N,s,q,\vec k)>0,  $$
with $\omega_0(N,s,q,\vec k)$ a constant independent of $r$.
Using the expansion (\ref{minenergy}) of the minimizing solution
in $\BPP(1,{\pi\over Hp},{\pi\over Hp},\ZZ)$
we obtain the first alternative in the statement of the Theorem.
When $q\in {\pi\over Hp}\ZZ$ the remaining statements
follow as a corollary to Theorem~\ref{Thm3}.
\QED

\begin{rem}\rm
The result contained in Theorems~\ref{bigthm} and \ref{Thm3},
namely that the smallest possible energy per unit
cross-sectional area is obtained with $N=1$, $s=q_1$, $m=1$,
confirms the prediction of a period-$2p$ in $z$ staggered lattice solution
made by  Bulaevski\u\i--Clem \cite{BC}.
In that paper the authors indeed claim that it minimizes
energy among competing configurations, but only one other
vortex lattice is treated in their paper (a  period-4 lattice),
and no indication is provided as to how they deduced the
geometry of their solution.
Here we have shown much more:
we know that it is the minimizer among {\it all}
periodic solutions, in the regime $r<<1$.

The
``vortex plane solution'' of Theorodakis--Kuplevaksky (\cite{Th}, \cite{K}) is
obtained by taking $\delta_n=0$ for all $n$ (when
the value of $s$ permits such a choice.)   Note that such a choice
{\it maximizes} $\Omega^{(2)}$ on
the constraint set $\SSS_r$.   Since these would also constitute
non-degenerate critical points of the finite-dimensional
function $g(r,\delta_2,\dots,\delta_N)$, by the above arguments
they describe {\it bona fida} solutions to the
  Lawrence--Doniach system with periodic boundary conditions,
but they are unstable.  
\end{rem}
\begin{rem}\rm  When $s\neq jq_1$, $j\in \ZZ$ the
trivial choice of minima $\delta_n-\delta_{n-1}=\pm \pi$
is not admissible, and the lattice is ``frustrated''.  
As previously mentioned, in this generic case 
 the minimizer of $F(N,s)$ might be non-unique, leading to different
asymptotics for the minimizing solutions along different subsequences
$r\to 0$.  Nevertheless, the dependence of the energy minimizers
 on subsequential limits $r\to 0$ could be eliminated 
when it is known that the absolute minimizer 
$(\delta_2^*,\dots,\delta_N^*)$ of $F(N,s)$ is unique and
non-degenerate. 
An example is when the minimizer of
the finite dimensional problem (\ref{FDprob})
$(\delta^*_2,\dots,\delta^*_N)$ satisfies
\be\label{allneg}   C_n:= \cos(\delta^*_n-\delta^*_{n+1})<0 \quad
    \mbox{for all $n=1,\dots,N$.}
\ee
This will be the case when the parallelogram $\Pi_{N,s,q}$
is very close to the optimal ones described by (\ref{even})
and (\ref{odd}).
Assuming (\ref{allneg}) holds for the
minimizer, a simple calculation shows that the Hessian 
is the ($(N-1)\times (N-1)$) tridiagonal, symmetric matrix 
$D^2 g(0,\delta^*_2,\dots,\delta^*_N)=: [M_{m,n}]_{m,n=2,\dots,N}$ with 
$$    M_{n,n+1}= C_n \quad n=2,\dots, N-1;
      \qquad M_{n,n}= -(C_{n-1} + C_n), \quad n=2,\dots,N. 
$$
A null vector $\vec v=(v_2,\dots, v_N)$ of $M$ satisfies:
\beann
&&  C_n (v_{n+1}-v_n)  - C_{n-1} (v_n - v_{n-1})=0, 
           \qquad n=3,\dots,N-1; \\
&&  C_2(v_3-v_2) - C_1 v_2 =0, \qquad
           C_N v_N - C_{N-1} (v_N- v_{N-1}) =0.
\eeann
If $v_2=0$ then clearly $\vec v=0$, so we may assume
$v_2>0$, in which case equations $n=2,\dots,N-1$ imply
$0<v_2<\dots<v_N$.  But this contradicts the $n=N$ equation,
and therefore the only solution is the trivial one.
In conclusion the minimizer is non-degenerate,
and we can repeat the same arguments as in the
cases (\ref{even}) and (\ref{odd}) to conclude uniqueness
for energy minimizers.
\end{rem}

\section{The periodic finite layer case}
\setcounter{thm}{0}

Finally, we consider the case of a finite number of planes $N$,
each of infinite extent in $x$ and $y$, assuming that the
currents and field strength are periodic functions in $x$.
Since this case is very similar to the doubly periodic case
treated in the previous sections we give an outline of how to modify
the formulation of the problem and its solution to fit this
somewhat simpler case.

In the finite-layer case the `t~Hooft condition is greatly simplified,
since (by the argument of Theorem~\ref{bipercoulomb})
we may take $\vec A$ to be periodic in $x$.  
Let $\vec k=(k_n)_{n\in\ZZ}$ with $k_0=0$, $\vec k\neq \vec 0$.
We say that $(f_n,\phi_n,\vec A)$ belongs to the periodic class 
${\cal  P}={\cal P}(q,\vec k)$
 if there exists a constant $\omega\in\RR$ such that:
$$
 \vec A\in H^1_{loc}(\RR^2;\RR^2), \  f_n\in H^1_{per}(\RR), \ 
        \phi_n\in H^1_{loc}(\RR); $$ 
\be
 \label{P2}
  \phi_n(x+2q)=\phi_n(x) + \omega + 2\pi k_n, \quad n=0,\dots, N;  
\ee
\be\label{P0} 
\int_{-q}^q\phi_0(x)\, dx=0=\int_{-q}^q\phi_1(x)\, dx;
\ee
\be
\label{P3} 
 \left. \begin{array}{c}
\vec A(x,z)=(H z,0) + (\xi_z,-\xi_x), \mbox{ with} \  \xi\in H^2_{loc}(\RR^2),\\ 
  \\
 \xi(x+2q,z)=\xi(x,z),\quad \xi(x,0)=\xi(x,Np)=0.
\end{array}  \right\}
\ee
If $\vec k=\vec 0$ we define ${\cal P}(q,\vec 0)$
by omitting (\ref{P0}).
Following the proof of Theorem~\ref{bipercoulomb},
 any configuration satisfying a single
`t~Hooft condition (in $x$ with $\ox$) as in (\ref{tHooft1}),
(\ref{tHooft2}) is gauge--equivalent to an element of ${\cal P}(q,\vec k)$.
Furthermore, in this class the $L^2$ norm of  $\curl\vec A$ 
controls $\vec A$ in $H^1$. 
 
\smallskip

As in the doubly periodic case, the observables are determined
entirely by the values of $(f_n,\phi_n,\vec A)$ in a single
period $x\in [-q,q]$, but now the number of planes is finite 
(indexed by $n=0,\dots,N$) and
each is independent.   Hence we define the Gibbs free energy  by
integration over a single period and summation over the $N+1$ planes,
 \beann
\Omp (f_n,\phi_n,\vec A)&=& \int_{-q}^{q} p\sum_{n=0}^N \left[
       \frac 12  (f_n^2-1)^2 + {1\over\kappa^2} (f'_n)^2
            + {1\over\kappa^2}(\phi'_n- A_x(x,z_n))^2 f_n^2
             \right]\, dx \\   
             &&\qquad + {r\over 2} \int_{-q}^q p\sum_{n=1}^N
                 \left( f_n^2 +f_{n-1}^2 
              -2 f_n f_{n-1}\cos(\Phi_{n,n-1} )\right) \, dx \\  
            &&\qquad  +\ {1\over  \kappa^2}
                \int_0^{Np}\int_{-q}^q \left(  {\partial A_x\over \partial z} 
                    - {\partial A_z\over \partial x}
      - H\right)^2 \, dx\, dz .
  \eeann 
By taking variation of $\Omp (f_n,\phi_n,\vec A)$ in the 
subclass of periodic functions, we see that the Euler--Lagrange equations
are exactly the same as for the finite width sample case,
see (2)--(5), (7) of \cite{ABB1}.
That is, they coincide with (\ref{bfeqn})--(\ref{cu})
for $n=1,\dots,N-1$, and the equations involving the top
and bottom surfaces $n=0,N$ must be modified to reflect
the fact that these two planes have only one ``nearest neighbor''.
\medskip

We may now continue as in the doubly periodic case.
By Lemma~\ref{inf_per} the minimum energy will be of
order $r$ if and only if we choose $q=m{\pi\over Hp}$
and $\vec k = m\ZZ$, with $m=1,2,\dots$. 
Furthermore, with these choices the minimizers
will have $f_n=|\psi_n|\sim 1$, and the use of
polar coordinated for the order parameter is well justified. We 
choose $q={\pi\over Hp}$, $\vec k=\ZZ$, and write
$\PP:={\cal P}({\pi\over Hp},\ZZ)$ in the following.

For the number $N$ of planes, $\kappa$, and $H$ fixed,
define the free energy
per unit area in a period strip of width $2q$ 
with winding numbers $\vec k$ by:
$$
   \epsilon_r(q,\vec k):=  {1\over 2qNp}\inf\left\{
       \Omp(f_n,\phi_n,\vec A): \ (f_n,\phi_n,\vec A)\in 
                            {\cal P}(q,\vec k)\right\}.
$$

\medskip

We obtain the following result:
\begin{thm}\label{Thm2}
Let $N,H,\kappa$ be fixed.
\begin{enumerate}
\item[(a)]   If there exists $m\in\NN$ such that
\be\label{goodchoice} q=q_m:= {m\pi\over Hp}
  \qquad\mbox{and}\qquad \vec k = m\ZZ,  
\ee
then there exists $r_0=r_0(N,H,\kappa,m)>0$ such that
for all $0<r<r_0$, $\epsilon_r(q,\vec k)=\epsilon_r({\pi\over Hp},\ZZ)$,
and the minimizers of $\Omp$ in ${\cal P}(q_m,m\ZZ)$
coincide with the  minimizers of $\Omp$ in $\PP$.
\item[(b)]  For any other choice of $q$, $\vec k$ there
exist constants $r_1,\omega_1>0$ (depending
on $N,H,\kappa,q,\vec k$) such that
$\epsilon_r(q,\vec k)\ge \omega_1$ for all $0<r<r_1$.
\end{enumerate}
\end{thm}

We note that, as in the bi-periodic case, we may obtain
an expansion of the minimizing solution in powers of $r$.
The minimizer in the special space $\PP$ (which gives
the absolute minimum of energy per unit area)
coincides with the period-$2p$ in $z$ lattice found
in the bi-periodic case at order $r$, except for an
edge effect in the order parameter in the top and
bottom planes.  More precisely, the fields and
currents $h(x,z)$, $j_x^{(n)}$, $j_z^{(n)}$ in a finite
stack still
 satisfy (\ref{per_2}) for each $n=0,\dots,N$,
and $f_n$ coincides with the expression given
in (\ref{per_2}) for $n=1,\dots, N-1$ but is modified
(see (\ref{edge}) and (\ref{edge0}) below) at
order $r$ for $n=0$ and $n=N$.

\medskip

The proof of Theorem~\ref{Thm2} follows
almost line-for-line the degenerate
perturbation procedure of the previous sections.
In particular,  the minimum
value of $\Omp$ at $r=0$ is attained by elements in the 
$(N-1)$-dimensional hyperplane
\bea\label{spmin0}
      \SSS & := &\{(f_n,\phi_n,\vec A)\in \PP: \ f_n\equiv 1, \  
      \phi_n(x)=\alpha_n +nHpx, \\
       & & \qquad   \  A_x=Hz, \ A_z=0, \ 
       \mbox{where $\alpha_0=\alpha_1=0$,\ \mbox{and}\
       $\alpha_2,\dots,\alpha_N\in{\bf R}$.} 
              \}
\nnn
\eea

The only significant difference with the doubly periodic
problem occurs at the top and bottom layers, $n=0,N$.  The
first-order expansion of the solution in $r$ requires that
these be treated differently.  For example
 the equation of current conservation (\ref{proj2}) now gives 
\be\label{p1}
v_{n,1}' - a_{x,1}(x,z_n) = 
\cases{
C_n - {\kappa^2\over 2Hp} 
   \cos(\delta_n + Hpx)
      + {\kappa^2\over 2Hp} \cos(\delta_{n+1} + Hpx), & $n=1,\dots N-1$; \cr
    C_0 + {\kappa^2\over 2Hp} \cos(\delta_1 + Hpx), & $n=0$; \cr
    C_N - {\kappa^2\over 2Hp} 
   \cos(\delta_N + Hpx), & $n=N$;
    }
\ee
The order $r$ term in magnetic field $b_1$ is exactly as in the doubly periodic
case and is given by (\ref{p2}) (with 
as yet undetermined constants $D_1,\dots,D_N$.)
Therefore, using equation (\ref{p2}) and the jump conditions (\ref{proj4}),
we obtain the following conditions on the constants $C_n$ and $D_n$:
\bea\label{p3}
&& D_{n+1} -2 D_n + D_{n-1} = p^2 D_n,\qquad n=1,\dots,N-1,  \\
\label{p3.1} && D_2-2D_1 = p^2 D_1, \quad D_{N-1}-2D_N =p^2 D_N, \\
\nnn && -pC_n=(D_{n+1}-D_n),\qquad n=1,\dots,N-1,  \\
\nnn && -pC_0=D_1,\quad pC_N=D_N.
\eea
The maximum principle for the second-order difference
equation (\ref{p3}) with boundary condition (\ref{p3.1})
implies that the unique solution is $D_n=0$,  $C_n=0$
for all $n$.
In consequence the perturbed manifold $\SSS_r$
consists of the  same configurations (\ref{persol})
as for the bi--periodic case, except for the superconducting order
parameter which  coincides with the expression given (\ref{persol}) for
$n=1,\dots,N-1$, but for $n=0$ or $n=N$ we obtain
\bea\label{edge}
f_N&=&1\, +\, {r\over 2}
\left(-\frac12 + {\kappa^2\over 2(H^2p^2 + 2\kappa^2)}
     \cos(\delta_N + Hp x) \right)
      \, + \, O(r^2),  \\
\label{edge0}
f_0&=&1\, +\, {r\over 2}
\left(-\frac12 + {\kappa^2\over 2(H^2p^2 + 2\kappa^2)}
     \cos(Hp x) \right)
      \, + \, O(r^2).
\eea

 Substitution of (\ref{persol}), (\ref{edge}), (\ref{edge0}) into 
$\Omp$ leads to an expansion of the energy in the same
form as (\ref{energy2}),
$$
\Omp(\sigma + w(r,\sigma)) = 2Npq\left\{ r 
   + r^2\left( \tilde C_0 + \tilde C_1\, {1\over N}\sum_{n=1}^{N-1}  
      \cos(\delta_n - \delta_{n+1}) 
           \right) \right\}
       + O(r^3) ,
$$
except for the constants $\tilde C_0\in \RR$ and $\tilde C_1>0$
which differ from the doubly periodic case due to the slightly
different form of solutions for the top and bottom planes.
The significant difference from the previous case is that
when there are finitely many planes there is only
 the single constraint $\delta_1=0$
 (which comes from removing the translation invariance in the
 definition of the space $\PP$).  By the same arguments
as in the previous section the minimizer of $\Omp$
will be determined by minimizing the leading term in
the energy expansion,
$$  G(\delta_2,\dots,\delta_N)=\sum_{n=1}^{N-1} \cos(\delta_n - \delta_{n+1}) ,
\quad
\delta_1=0. 
$$ 
By inspection, the minimizer
 is obtained by choosing $\delta_{n+1}-\delta_n=\pm \pi$;
 for example $\delta_1=\pi$, $\delta_2=0$, $\delta_3=\pi$, \dots
Unlike the doubly periodic case we are always free to make this
choice, and it is easy to verify that this configuration gives
a non-degenerate minimizer of the finite-dimensional
function $g(0,\sigma)$, $\sigma=(\delta_2,\dots, \delta_N)$.
In particular
the period-$2p$ in $z$, period-$2q_1$ in $x$
lattice appears naturally (at order $r$) in the solutions.
(See figure~1.)  As in the finite-width case the solutions
feel the top and bottom edges only in the first  plane and
first gap at order $r$; expansion to higher orders will
reveal the effect of the finite stack at order $r^k$ in the
$k$th stack from the top or bottom.

\section{Conclusions}

Finally, we summarize our results on the periodic problem, and
compare with the (quite different) conclusions obtained
for finite-width samples in our previous paper \cite{ABB1}.

As we have seen, in the limit $r\to 0$ energy minimization
selects a preferred period geometry and quantized flux from
all possible periodic configurations represented by the
spaces $\BPP(N,s,q,\vec k)$.  The optimal solution is
${2\pi\over Hp}$-periodic in $x$, and repeats itself with
a horizontal shift of a half-period in $x$ when we climb
from one superconducting plane to the next.  Each period
parallelogram contains exactly one quantum $2\pi$ of flux,
and one Josephson vortex per period in each gap between
the superconducting planes.
 The Josephson
lattice geometry is the same as was predicted by Bulaevsk\u \i i
\& Clem \cite{BC}, but with our rigorous analytical approach
we may now assert that it is the unique energy minimizing periodic
configuration (among all possible geometries)
for all sufficiently small $r$.  

In the finite-width case \cite{ABB1} the conclusions were surprisingly
very different.  For {\it any} finite width sample,
$-L\le x\le L$, when $\sin (HpL)\neq 0$ the unique
energy minimizer for $r\sim 0$ is a vortex plane configuration, with
Josephson vortices vertically aligned, and the magnetic
field approximately uniform in $z$ (except for edge effects
at the top and bottom of a sample of finitely many planes.)
The exceptional values of the applied field $H$ for which
$\sin(HpL)=0$ correspond to first-order phase transitions
occuring when a new vortex plane is nucleated into the
sample from the lateral edges.

We can explain this apparent conflict by
examination of the expansion of the energy
in powers of $r$ near the degenerate manifold $\SSS$.
In (23) of \cite{ABB1} vortex planes are preferred
at order $r$ in the energy because of a {\it surface} term,
 a quantity which scales like the length $Np$ of
the lateral edges of the sample cross-section.  In the periodic case this
term does not appear, and in (\ref{energy2}) the
distinction between lattice geometries appears at order
$r^2$, in a term which scales as the cross-sectional area
$2qNp$ of the sample.  For any size sample, by making
$r$ small enough the order $r$ surface term will dominate,
but the value of $r$ must decrease if the surface term is
to continue to prevail with increasing sample width $L$.
If we try to keep $r$ fixed while increasing $L$ then inevitably
the order $r^2$ term will compete with the order $r$ term.
At that point the perturbation expansion will surely have
lost its validity.

This general argument is supported by
the analysis of the range of validity of the expansion
in Lemma~\ref{Lya}, presented in section~5 of \cite{ABB1}.
Indeed, it seems clear that the radius
of validity of the expansions in $r$  deteriorates with increasing
sample width $L$. (See Remark~5.3 of \cite{ABB1}.)

If we approximate a macroscopic sample by an infinite one and
seek periodic solutions, the period plays the role of the sample
width in the estimates of the interval of validity of section ~5
in \cite{ABB1}.  Since the period of the absolute minimizer
is given by ${2\pi\over Hp}$ we can expect the interval of
validity to extend to physically appropriate values of $r$
when the applied field $H$ is large enough.  Therefore we
may apply our analysis to the transparent state of the 
high-$T_c$ superconductors in high external fields.

\begin{figure}
\centering
\includegraphics{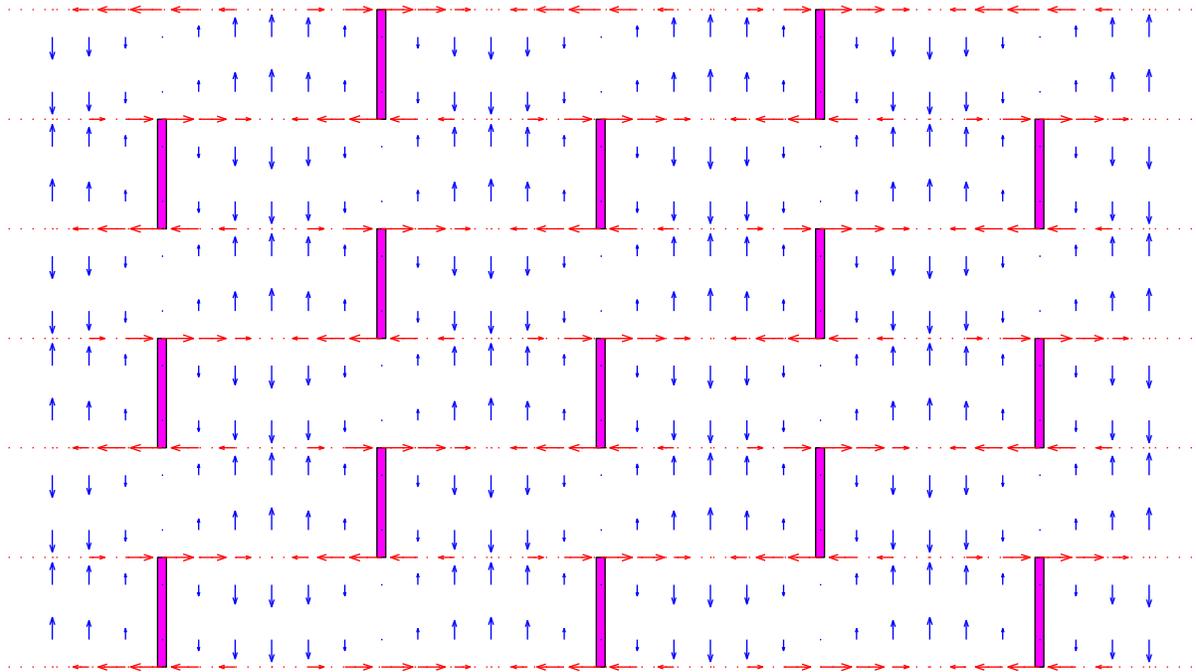}
\caption{{\small \sl
Period-2$p$ (in $z$) vortex lattice, for a sample with a 
finite number of superconducting planes (Indicated
by horizontal dotted lines.)
Horizontal arrows indicate the
in-plane currents
$j_x^{(n)}$ and the vertical arrows depict the Josephson
currents $j_z^{(n)}$ between adjacent planes.
The magnetic field
$h(x,z)$ and supercurrents $j_x$, $j_z$ are periodic with
period $2p$ in the $z$-direction
(and period ${2\pi\over Hp}$ in
$x$).  The vortices (local maxima of
$h$) lie along the starred segments, and form a staggered lattice.
If we choose the midpoint of each segment to label each vortex
the resulting lattice is diamond-shaped.}}
\end{figure}

\end{document}